\documentclass[11pt]{article}

\usepackage[utf8]{inputenc}
\usepackage[T1]{fontenc}
\usepackage[USenglish]{babel}
\usepackage[a4paper,margin=32mm,headheight=14pt]{geometry}
\usepackage{lmodern}
\usepackage{microtype}
\usepackage{csquotes}
\usepackage{amsmath,amssymb}
\usepackage{amsthm}
\usepackage{mathtools}
\usepackage{booktabs}
\usepackage{needspace}
\usepackage{xcolor}
\usepackage{tikz}
\usepackage{placeins}
\usetikzlibrary{calc,patterns,arrows.meta,decorations.pathreplacing}
\usepackage{fancyhdr}
\usepackage[
  colorlinks=true,
  linkcolor=blue!45!black,
  citecolor=green!35!black,
  urlcolor=blue!55!black,
  pdfauthor={Alex Junior Gomez Saltachin},
  pdfcreator={LaTeX with hyperref},
  pdftitle={Conjecture B for X10 and its genus-two mirror},
  pdfsubject={Algebraic geometry and mirror symmetry},
  pdfkeywords={Orbifold del Pezzo surfaces, maximally mutable Laurent polynomials, mirror symmetry, quantum periods, X10}
]{hyperref}
\usepackage[capitalize,nameinlink]{cleveref}

\theoremstyle{plain}
\newtheorem{theorem}{Theorem}[section]
\newtheorem{proposition}[theorem]{Proposition}
\newtheorem{lemma}[theorem]{Lemma}
\newtheorem*{conjectureB}{Conjecture B}
\theoremstyle{definition}
\newtheorem{definition}[theorem]{Definition}
\theoremstyle{remark}
\newtheorem{remark}[theorem]{Remark}

\crefname{theorem}{Theorem}{Theorems}
\crefname{proposition}{Proposition}{Propositions}
\crefname{lemma}{Lemma}{Lemmas}
\crefname{definition}{Definition}{Definitions}

\newcommand{\cX}{\mathcal X}
\newcommand{\CT}{\operatorname{CT}}
\newcommand{\Eff}{\operatorname{Eff}}
\newcommand{\A}{\mathbb A}
\newcommand{\C}{\mathbb C}
\newcommand{\Q}{\mathbb Q}
\newcommand{\Z}{\mathbb Z}
\newcommand{\Pj}{\mathbb P}
\newcommand{\1}{\mathbf 1}

\usepackage[
  backend=biber,
  style=numeric,
  sorting=nyt,
  giveninits=true,
  maxbibnames=99,
  doi=true,
  url=false,
  isbn=false
]{biblatex}
\DeclareRobustCommand{\wangmirrorcite}{\cite[Theorem 1.1]{Wang}}

\fancypagestyle{plain}{%
  \fancyhf{}
  \fancyfoot[C]{\thepage}
  
  }
\begin{document}

\title{\textbf{Conjecture B for \texorpdfstring{$X_{10}$}{X10} and its genus-two mirror}}
\author{Alex Junior Gomez Saltachin\\[0.5em]
\small Departamento de Matem\'atica, Pontif\'icia Universidade Cat\'olica do Rio de Janeiro\\
\small Rua Marqu\^es de S\~ao Vicente, 225, Rio de Janeiro, RJ 22451-900, Brazil\\
\small \href{mailto:alex.gomez@pucp.edu.pe}{alex.gomez@pucp.edu.pe}}
\date{}

\maketitle

\begin{abstract}
We prove Conjecture B of Akhtar--Coates--Corti--Heuberger--Kasprzyk--Oneto--Petracci--Prince--Tveiten unconditionally for the normalized family of maximally mutable Laurent polynomials with prescribed $T$-binomial edge coefficients associated with the degree-$10$ orbifold del Pezzo hypersurface $X_{10}\subset\Pj(1,2,3,5)$.  An explicit $S$-extended quantum-period calculation determines the affine mirror map $\lambda=a-346$ and proves that the regularized quantum period equals the classical period as a formal power series throughout the family.
\end{abstract}



\section{Introduction}

The mirror-symmetry program for orbifold del Pezzo surfaces associates mutation classes of Fano polygons with $\mathbb Q$-Gorenstein ($qG$) deformation families of del Pezzo surfaces of class TG, meaning del Pezzo surfaces with cyclic quotient singularities that admit a $qG$-degeneration with reduced fibers to a normal toric del Pezzo surface \cite[Definition 1]{ACC+}.  Conjecture~B predicts a period-preserving affine correspondence between families of maximally mutable Laurent polynomials (MMLPs) and the twisted parameter space of a generic $qG$-deformation.

This paper proves that prediction for the general degree-$10$ hypersurface $X_{10}\subset\Pj(1,2,3,5)$.  We write $\cX_{10}$ for its canonical smooth Deligne--Mumford stack.  The surface has a unique singular point of type $\frac13(1,1)$.  Its mirror mutation class contains a polygon that Oneto--Petracci call n.1 \cite[entry 11]{OP}, namely $P=\operatorname{conv}\{(-3,-5),(-3,5),(7,5)\}$.
The generic $qG$-deformation of the toric surface associated with $P$ belongs to the family $X_{1,1/3}$, identified by Corti--Heuberger \cite[Section 7 and Table 4, entry 1]{CortiHeuberger}.  Oneto--Petracci represent its general member by the hypersurface above \cite[entry 11]{OP}.  The polygon has singularity content $(10,\{\frac13(1,1)\})$ \cite[Table 4 and Theorem 9]{KNP}.  These identifications place the hypersurface and its mirror polygon in the setting of Conjecture~B.

\Needspace{10\baselineskip}
Oneto--Petracci studied the specialization $a_{[4,3]}=360$ of their non-extended twisted $I$-function.  Its Laurent polynomial is $f_{360}(x,y)=\frac{y^5}{x^3}\left(1+x+\frac1y\right)^{10}-2520$, and its regularized quantum period begins $1+32448360t^2+515050578000t^3+\cdots$ \cite[entry 11]{OP}.  They note earlier calculations of this specialization by Kasprzyk--Tveiten \cite[Section 2]{OP}.  Their Theorem 2.2 proves the period-preserving conclusion under the conjectural orbifold extensions of quantum Lefschetz (Conjecture 5.2) and the Abelian/non-Abelian correspondence (Conjecture 6.1) \cite[Theorem 2.2 and Conjectures 5.2 and 6.1]{OP}.  For the weighted-hypersurface calculation here, the relevant conjectural input is Conjecture 5.2.  We establish the comparison using Wang's mirror theorem without that assumption.

Sultani--Webb give a separate rigorous computation of the full quantum period for $X_{1,7/3}$, another surface in the Oneto--Petracci classification \cite[Theorem 1.1.3]{SultaniWebb}.  Their nonabelian CI-GIT presentation is related to a weighted Grassmannian and uses quasimap $I$-functions and wall-crossing.  For $X_{10}$, we use a toric $S$-extension with Wang's mirror theorem, which applies without the convexity assumption.

This case is especially tractable: the polygon has a single residual point, while the canonical stack has a unique sector of age less than one.  Both parameter spaces are therefore one-dimensional.  The $S$-extended GIT presentation makes the scalar quantum period explicit, and the constant-term comparison yields the period-preserving affine correspondence.

Here $L_P^T$ denotes the affine space of MMLPs with Newton polygon $P$, normalized vertex coefficients, zero constant term, and the prescribed $T$-binomial edge coefficients.  We use $a=[x^4y^3]f$ as its affine coordinate.

\begin{theorem}[Completeness of $L_P^T$]\label{thm:mmlp-family}
Set $H=1+x+y^{-1}$.  The family $L_P^T$ is
\begin{equation}\label{eq:fa}
\boxed{
 f_a(x,y)=
 x^{-3}y^5H^{10}-2520
 +(a-360)\left(x^{-2}y^3H^6-60\right),
 \qquad a\in\A^1.
}
\end{equation}
Moreover, the coefficient of $x^4y^3$ in $f_a$ is $a$.
\end{theorem}

Center this coordinate by setting
\[s:=a-360,\qquad C(s):=2520+60s.\]
With $A=x^{-3}y^5H^{10}$ and $B=x^{-2}y^3H^6$, the same family member is
\[
f_s:=f_{s+360}=A+sB-C(s).
\]
The quantum computation in Section~\ref{sec:Ifunction} will determine the twisted coordinate $\lambda=s+14$.  We use $f_a$ when the distinguished coefficient is relevant and $f_s$ for computations in the centered coordinate.

The polygon $P$ has $36$ interior lattice points, so a sufficiently nondegenerate Laurent polynomial with Newton polygon $P$ has a general fiber of genus $36$.  The MMLP mutability conditions force boundary degeneracies instead.  Coates--Kasprzyk--Pitton--Tveiten (CKPT) show that the genus of a general MMLP fiber is one plus the number of residual points in $P^\circ$ \cite[Theorem 3.13]{CKPT}, which gives genus two here.  The common factor $H^6$ in $A$ and $B$ makes this degeneration explicit and motivates the substitutions in the hyperelliptic construction.  The factorized family and its genus-two model were obtained earlier in the author's categorical-genus preprint \cite[Proposition 7.1]{GomezCatGenus}.

\Needspace{12\baselineskip}
The following theorem records that model in the present notation.
\begin{theorem}[Explicit genus-two model]\label{thm:hyperelliptic}
For $s,\eta\in\A^1$, consider the fiber $f_s(x,y)=\eta$ and set $T=\eta+C(s)$.  The generic fiber is birational to
\begin{equation}\label{eq:hyperelliptic}
\boxed{
 v^2=T(u+s)\left(4u^3(u+s)^2-Tu^2(u+s)+4T^2\right).
}
\end{equation}
For general $(s,T)$ the polynomial on the right is squarefree of degree six, so the smooth projective model of the fiber has genus two.
\end{theorem}

To state the period correspondence, recall that $\cX_{10}$ is the canonical smooth Deligne--Mumford stack defined above, and denote its unique sector of age $\frac23$ by $\alpha$.  Its twisted parameter space is $H^{\mathrm{ts}}_{\cX_{10}}=\C\1_\alpha$.  Write $\widehat G_{\cX_{10}}$ for the regularized quantum period and $\pi_f$ for the classical period of a Laurent polynomial $f$.  Section~\ref{sec:background} gives the precise definitions of these periods and the twisted parameter space.  The quantum calculation applies Wang's mirror theorem to the $S$-extended GIT presentation to determine the affine correspondence between the parameter spaces and prove the following period identity.

\begin{theorem}[Conjecture B for $X_{1,1/3}$]\label{thm:conjB}
Let $f_a$ be the family \eqref{eq:fa}.  The affine map
\[
\Phi:L_P^T\longrightarrow H^{\mathrm{ts}}_{\cX_{10}},
\qquad
\Phi(f_a)=(a-346)\1_\alpha,
\]
is an affine isomorphism and satisfies
\[
\widehat G_{\cX_{10}}\bigl(\Phi(f_a);t\bigr)=\pi_{f_a}(t).\]
Consequently, Conjecture B of \cite{ACC+} holds for polygon n.1 and the family $X_{1,1/3}$.
\end{theorem}

Theorem~\ref{thm:mmlp-family} follows from Lemmas~\ref{lem:edge-mutability} and \ref{lem:direct-T-space}.  Section~\ref{sec:mmlp} derives the birational model in Theorem~\ref{thm:hyperelliptic}.  Theorem~\ref{thm:conjB} follows from the applicability result in Lemma~\ref{lem:Wang-hypotheses}, the effective-degree and leading-exponent analyses in Lemmas~\ref{lem:effective-semigroup} and \ref{lem:z-order}, the scalar completeness result in Proposition~\ref{prop:negative-ell}, and the final multinomial comparison.

\section{Background and conventions}\label{sec:background}

We recall the conventions needed to formulate Conjecture B.  On the classical side, these concern Fano polygons and maximally mutable Laurent polynomials with prescribed $T$-binomial edge coefficients.  On the quantum side, they concern twisted sectors and orbifold quantum periods.

\begin{definition}[Fano polygons and singularity content]\label{def:fano-polygon}
Let $N\simeq\Z^2$ be a lattice with dual lattice $M$.  A \emph{Fano polygon} is a convex lattice polygon $P\subset N_{\mathbb R}$ whose vertices are primitive lattice points and whose strict interior contains the origin \cite[Definition 1.1]{CKPT}.  For an edge $E$ of $P$, let $v_E\in N$ be a primitive vector parallel to $E$ and $w_E\in M$ its primitive inward-pointing normal.  The \emph{local index}, or \emph{height}, $r_E$ is determined by
\[
\langle w_E,p\rangle=-r_E
\qquad (p\in E),
\]
and the \emph{lattice length}, or \emph{width}, of $E$ is the number of primitive lattice segments into which $E$ decomposes.  These conventions agree with the singularity-content formalism of \cite[Section 3(a)]{CKPT}.  Write
\begin{equation}\label{eq:width-decomposition}
\operatorname{width}(E)=n_Er_E+w_{0,E},
\qquad
0\le w_{0,E}<r_E.
\end{equation}
The cone $\operatorname{conv}\{0,E\}$ decomposes crepantly into $n_E$ \emph{primitive $T$-cones}, each of width and local index $r_E$, together, when $w_{0,E}$ is positive, with one \emph{residual cone} of width $w_{0,E}$ and local index $r_E$.  The \emph{singularity content} of $P$ consists of the total number $\sum_E n_E$ of primitive $T$-cones and the basket of residual cones \cite[Section 3(a), especially Definition 3.2]{CKPT}.  A choice of residual points is the set of lattice points attached to these residual cones in \cite[Definitions 3.5--3.6]{CKPT}.
\end{definition}

Let $X_P$ denote the toric surface whose fan is generated by the faces of $P$.  We use this edge decomposition below.

\begin{definition}[Normalized Laurent polynomials and the $T$-binomial parameter space]\label{def:mutability-T-space}
For a Laurent polynomial $f$, write $\CT(f)$ for its constant term and $[x^iy^j]f$ for the coefficient of $x^iy^j$.  A Laurent polynomial is \emph{normalized} if the coefficient at every vertex of its Newton polygon is one \cite[Definition 2.2]{CKPT}.  For an edge $E$, write
\[
f=\sum_h f_h,
\]
where $f_h$ is supported on the affine slice $\langle w_E,\cdot\rangle=h$.  Set
\[
a_E=(1+x^{v_E})^{n_E}.
\]
The CKPT edge-mutability condition with factor $a_E$ is the slice divisibility condition
\begin{equation}\label{eq:edge-mutability-condition}
a_E^{-h}\mid f_h
\qquad\text{for every }h<0,
\end{equation}
after identifying the coordinate ring of each slice with a monomial multiple of $\C[x^{\pm v_E}]$.  Equivalently, $f$ is mutable with respect to the edge mutation factor $a_E$ in the sense of \cite[Definition 1.7 and Proposition 3.7]{CKPT}.

Following Conjecture B of \cite{ACC+}, let $L_P^T$ be the affine space of normalized MMLPs with Newton polygon $P$, zero constant term, and prescribed $T$-binomial edge coefficients.  If $\tau=x^{v_E}$ is a monomial coordinate along an edge $E$ of height $r_E$ and width $n_Er_E+w_{0,E}$, its edge polynomial, after removing an endpoint monomial, has the coefficients of $(1+\tau)^{n_Er_E}$ when $w_{0,E}=0$ and of $(1+\tau)^{n_Er_E}(1+\tau^{w_{0,E}})$ otherwise.  These prescribed edge coefficients are separate conditions from vertex normalization.
\end{definition}

The edge-mutability factors account for powers such as $(1+x)^{2j}$ and $(1+y)^{-3i}$ that occur later.  Only the following consequences of the general MMLP theory will be used.

\begin{proposition}[MMLP facts]\label{prop:mmlp-facts}
Let $P$ be a Fano polygon.
\begin{enumerate}
\item A Laurent polynomial in two variables with Fano Newton polygon is maximally mutable if and only if it is a specialization of the general Laurent polynomial produced by Proposition 3.7 using its prescribed factors $a_E$, as characterized in Proposition 3.9 \cite[Propositions 3.7 and 3.9]{CKPT}.
\item If $f$ is an MMLP in two variables, then a general member of its associated pencil of curves has genus one plus the number of residual points in $P^\circ$ \cite[Theorem 3.13]{CKPT}.
\end{enumerate}
\end{proposition}

On the quantum side, let $\cX$ be a smooth Deligne--Mumford stack.

\begin{definition}[Inertia sectors, age, and twisted parameters]\label{def:twisted-parameters}
If a connected component $\gamma$ of its inertia stack acts on a tangent fiber with eigenvalues
\[
\exp(2\pi\sqrt{-1}\,\theta_1),\ldots,
\exp(2\pi\sqrt{-1}\,\theta_m),
\qquad 0\le\theta_i<1,
\]
then its \emph{age} is $\operatorname{age}(\gamma)=\sum_i\theta_i$.  We refer to nontrivial sectors whose age lies strictly between zero and one as sectors of age less than one.  We use the twisted parameter space
\[
H^{\mathrm{ts}}_{\cX}
=\bigoplus_{0<\operatorname{age}(\gamma)<1}\C\1_\gamma.
\]
For $\cX_{10}$, the nonidentity elements at its unique $\frac13(1,1)$ point have ages $2/3$ and $4/3$.  Thus the unique sector of age less than one, denoted by $\alpha$, spans the twisted parameter space $H^{\mathrm{ts}}_{\cX_{10}}=\C\1_\alpha$, which we identify with $\C$ by $\lambda\1_\alpha\mapsto\lambda$.
\end{definition}

The age grading and orbifold Gromov--Witten conventions are those of Chen--Ruan and Abramovich--Graber--Vistoli \cite{ChenRuan,AGV}.

\begin{definition}[Orbifold quantum period and regularization]\label{def:orbifold-period}
Following the convention of \cite[Definition 3.2]{OP}, for $\lambda\in H^{\mathrm{ts}}_{\cX}$ the orbifold quantum period is the genus-zero descendant generating series
\[
G_{\cX}(\lambda;t)=\sum_{r\ge0}c_r(\lambda)t^r,
\]
If a term is indexed by a stable-map degree $d$ and marked insertions $b_1,\ldots,b_n$ from sectors of age less than one, then its exponent of $t$ is
\begin{equation}\label{eq:period-degree}
\delta(d;b_1,\ldots,b_n)
=-K_X\mathbin{\cdot}d
+\sum_{i=1}^n\bigl(1-\operatorname{age}(b_i)\bigr).
\end{equation}
For $r\ge2$, Definition~3.2 of \cite{OP} gives
\[
c_r(\lambda)=
\sum_{\substack{d\in\Eff(\cX_{10}),\ n\ge0\\
\delta(d;\underbrace{\alpha,\ldots,\alpha}_{n})=r}}
\frac{\lambda^n}{n!}
\left\langle
\underbrace{\1_\alpha,\ldots,\1_\alpha}_{n},
\phi_{\mathrm{vol}}\psi^{r-2}
\right\rangle^{\cX_{10}}_{0,n+1,d},
\]
and $c_0(\lambda)=1$, $c_1(\lambda)=0$.  Here $\Eff(\cX_{10})$ is the semigroup of effective stable-map curve classes, $\phi_{\mathrm{vol}}$ is the point class, and $\psi$ is the cotangent-line class at the last marking.  The \emph{regularized quantum period} is
\[
\widehat G_{\cX}(\lambda;t)
=\sum_{r\ge0}r!\,c_r(\lambda)t^r.
\]
\end{definition}

\begin{definition}[Classical Laurent-polynomial period]\label{def:classical-period}
For a Laurent polynomial $f$, its classical period in the standard Fano normalization \cite{CCGK} is
\[
\pi_f(t)=\sum_{n\ge0}\CT(f^n)t^n.
\]
\end{definition}

With these conventions, Conjecture B takes the following form.

\begin{conjectureB}[\cite{ACC+,OP}]
Let $P$ be a Fano polygon and let $X$ be a generic $qG$-deformation of the toric surface $X_P$, with canonical smooth Deligne--Mumford stack $\cX$.  There is an affine-linear isomorphism
\[
\Phi:L_P^T\longrightarrow H^{\mathrm{ts}}_{\cX}
\]
such that, for every $f\in L_P^T$, the following equality holds:
\[
\widehat G_{\cX}\bigl(\Phi(f);t\bigr)=\pi_f(t).
\]
\end{conjectureB}

\section{The Laurent-polynomial mirror and its genus-two fibers}\label{sec:mmlp}\label{sec:genus2}

We derive $L_P^T$ by imposing the $T$-binomial slice conditions along the three edges of $P$.  The same factorization that makes these conditions explicit will yield the genus-two model of a general fiber.  The polygon is
\[
P=\operatorname{conv}\{(-3,-5),(-3,5),(7,5)\}.
\]
The left edge $x=-3$, top edge $y=5$, and diagonal edge $y-x=-2$ have lattice length $10$ and local indices $3$, $5$, and $2$, respectively.  Thus their widths decompose as
\[
10=3\cdot3+1,\qquad 10=2\cdot5,\qquad 10=5\cdot2.\]
Consequently the polygon has ten primitive $T$-cones and one residual cone of type $\frac13(1,1)$, in agreement with \cite{KNP}.

We introduce the Laurent polynomials
\[
H=1+x+y^{-1},\qquad
A=x^{-3}y^5H^{10},\qquad
B=x^{-2}y^3H^6.
\]
Their constant terms satisfy
\begin{equation}\label{eq:basic-CT}
\CT(A)=\frac{10!}{3!5!2!}=2520,
\qquad
\CT(B)=\frac{6!}{2!3!1!}=60.
\end{equation}
Their distinguished coefficients satisfy
\begin{equation}\label{eq:basic-a}
[x^4y^3]A=360,
\qquad
[x^4y^3]B=1.
\end{equation}
Accordingly, consider the candidate
\[
f_s=A+sB-C(s).\]
It has zero constant term, and its coefficient of $x^4y^3$ is $s+360$, which equals $a$.
At the vertices $(-3,5)$, $(-3,-5)$, and $(7,5)$, the coefficients of $A=x^{-3}y^5(1+x+y^{-1})^{10}$ are obtained by choosing, respectively, ten copies of $1$, ten copies of $y^{-1}$, and ten copies of $x$, so each coefficient is $1$.  The Newton polygon of $B$ is $\operatorname{conv}\{(-2,3),(-2,-3),(4,3)\}$.  None of the vertices of $P$ lies in this triangle, and the constant correction changes no vertex coefficient.

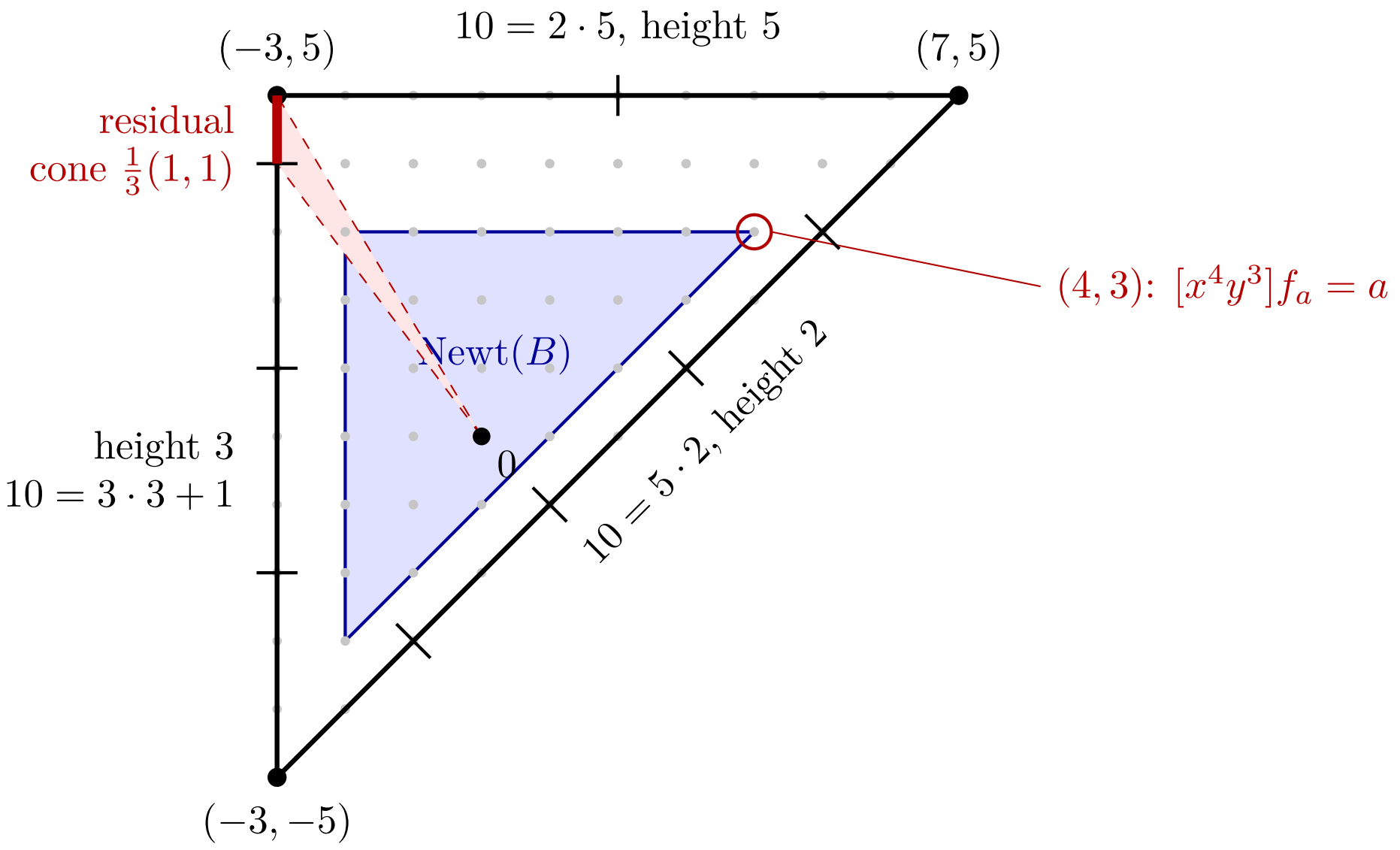
\begin{figure}[!htbp]
\centering
\begin{tikzpicture}[scale=0.60]
  % Newton polygon of B (shaded)
  \fill[blue!12] (-2,3) -- (-2,-3) -- (4,3) -- cycle;
  \draw[blue!60!black,thick] (-2,3) -- (-2,-3) -- (4,3) -- cycle;
  \node[blue!60!black,font=\small] at (0.2,1.2) {$\mathrm{Newt}(B)$};
  % residual cone at the origin
  \fill[red!10] (0,0) -- (-3,4) -- (-3,5) -- cycle;
  \draw[red!70!black,dashed] (0,0) -- (-3,4);
  \draw[red!70!black,dashed] (0,0) -- (-3,5);
  % lattice points of P
  \foreach \x in {-3,...,7}
    \foreach \y in {-5,...,5}{
      \pgfmathtruncatemacro{\d}{\y-\x}
      \ifnum\d>-3 \fill[gray!45] (\x,\y) circle (0.07); \fi
    }
  % polygon P
  \draw[very thick] (-3,-5) -- (-3,5) -- (7,5) -- cycle;
  % origin
  \fill (0,0) circle (0.13);
  \node[below right,font=\small] at (0,0) {$0$};
  % vertices
  \foreach \p in {(-3,-5),(-3,5),(7,5)} \fill \p circle (0.14);
  \node[above,font=\small] at (-3,5.15) {$(-3,5)$};
  \node[above,font=\small] at (7,5.15) {$(7,5)$};
  \node[below,font=\small] at (-3,-5.15) {$(-3,-5)$};
  % distinguished coefficient (4,3)
  \draw[red!70!black,thick] (4,3) circle (0.25);
  \draw[red!70!black] (4.25,3) -- (8.2,2.2) node[right,font=\small] {$(4,3)$: $[x^4y^3]f_a=a$};
  % top edge: height 5, 10 = 2*5
  \draw[thick] (2,4.7) -- (2,5.3);
  \node[above,font=\small] at (2,5.5) {$10=2\cdot 5$,\ height $5$};
  % diagonal edge: height 2, 10 = 5*2
  \foreach \px/\py in {-1/-3,1/-1,3/1,5/3} \draw[thick] (\px-0.25,\py+0.25) -- (\px+0.25,\py-0.25);
  \node[font=\small,rotate=45] at (3.3,-0.1) {$10=5\cdot 2$,\ height $2$};
  % left edge: height 3, 10 = 3*3+1, residual segment [4,5]
  \foreach \y in {-2,1,4} \draw[thick] (-3.3,\y) -- (-2.7,\y);
  \draw[red!70!black,line width=2.4pt] (-3,4) -- (-3,5);
  \node[left,font=\small,align=right] at (-3.4,-0.5) {height $3$\\ $10=3\cdot3+1$};
  \node[left,red!70!black,font=\small,align=right] at (-3.4,4.2) {residual\\ cone $\tfrac13(1,1)$};
\end{tikzpicture}
\caption{The polygon $P=\operatorname{conv}\{(-3,-5),(-3,5),(7,5)\}$ with its $66$ lattice points. Tick marks divide each edge into primitive $T$-cones: two of width $5$ on the top edge, five of width $2$ on the diagonal edge, three of width $3$ on the left edge, followed by the residual segment of width $1$ (red), which spans the cone of type $\tfrac13(1,1)$. The shaded triangle is $\operatorname{Newt}(B)$, which avoids all three vertices, and the circled point $(4,3)$ carries the distinguished coefficient.}
\label{fig:polygon}
\end{figure}

The binomial factorizations of $A$ and $B$ also make the $T$-binomial conditions visible on slices parallel to the three edges, which exhaust the boundary conditions for $P$.  Figure~\ref{fig:polygon} summarizes the polygon, its edge decomposition, and the distinguished coefficient.

\FloatBarrier

\begin{lemma}[Edge mutability]\label{lem:edge-mutability}
The family $f_s$ satisfies the $T$-binomial mutability conditions along all three edges of $P$.
\end{lemma}

\begin{proof}
For the top edge, the coefficient of $y^j$ in $A$ is obtained by choosing $5-j$ copies of $y^{-1}$ from $H^{10}$.  The remaining $5+j$ factors contribute $(1+x)^{5+j}$.  Thus the coefficient is given by
\[
[y^j]A=\binom{10}{5-j}x^{-3}(1+x)^{5+j}.
\]
The corresponding slice of $B$ is
\[
[y^j]B=\binom{6}{3-j}x^{-2}(1+x)^{3+j}.
\]
The inequalities
\[
5+j\ge2j\quad(j\le5),
\qquad
3+j\ge2j\quad(j\le3)
\]
give the required divisibility by $(1+x)^{2j}$ for $1\le j\le5$ and $1\le j\le3$, respectively.  The exponent $2j$ records the two primitive $T$-cones on this edge.

Along the left edge, choosing $i+3$ copies of $x$ in $A$ leaves $7-i$ copies of $1+y^{-1}$, while the corresponding exponent for $B$ is $4-i$.  Up to monomial factors, the negative slices $i=-3,-2,-1$ of $A$ therefore contain
\[
(1+y)^{10},\quad(1+y)^9,\quad(1+y)^8,
\]
and hence are divisible by the required powers $(1+y)^9,(1+y)^6,(1+y)^3$.  The two negative slices of $B$ contain $(1+y)^6$ and $(1+y)^5$.  Since a slice of grade $i<0$ must be divisible by $(1+y)^{-3i}$, these satisfy the same condition.

For the diagonal edge, we put $\xi=xy$, which gives $x=\xi/y$.  Under this substitution, the polynomials become
\[
A=\xi^{-3}y^{-2}(1+\xi+y)^{10},
\qquad
B=\xi^{-2}y^{-1}(1+\xi+y)^6.
\]
The edge lies on $j-i=-2$, so the $y$-exponent is the grading normal to it.  The grade $-2$ term is proportional to $(1+\xi)^{10}$, while at grade $-1$ the contributions from $A$ and $B$ are proportional to $(1+\xi)^9$ and $(1+\xi)^6$.  These contain the powers required by the five primitive $T$-cones on the diagonal edge.
The calculations verify the required divisibility for $A$ and $B$ separately along all three edges.  The constant polynomial $1$ has no negative-grade slices, so $A$, $B$, and $1$ each satisfy all three $T$-binomial conditions.  Since the conditions are linear and $f_s$ is the linear combination $A+sB-C(s)$, the same is true of $f_s$.
\end{proof}

Let $V$ be the vector space of Laurent polynomials supported on $P$ that satisfy the three $T$-binomial conditions.

\begin{lemma}[Direct description of the $T$-binomial space]\label{lem:direct-T-space}
For polygon n.1, the space $V$ has the following description:
\[
V=\operatorname{span}\{A,B,1\}.
\]
\end{lemma}

\begin{proof}
We use that a Laurent polynomial divisible by $(1+t)^e$ and supported in at most $e$ consecutive exponents is zero.  If the allowed exponents are exactly $e+1$ consecutive integers, divisibility forces it to be a scalar monomial multiple of $(1+t)^e$.

The polygon is described by
\[
x\ge-3,\qquad y\le5,\qquad y-x\ge-2.
\]
For a diagonal slice $y-x=q$, use $\xi=xy$ as its monomial coordinate.  The $T$-binomial exponents on the top, left, and diagonal edge slices are respectively $2j$, $-3i$, and $-5q$.  Starting with the top slice and then removing each slice that vanishes gives the following surviving parameter ranges.  The counts are read from the displayed inequalities after deleting points on earlier slices.
\[
\begin{array}{c|c|c|c|c}
\text{slice} & \text{required factor} & \text{exponents} & \text{count} & \text{consequence}\\
\hline
y=5 & (1+x)^{10} & x=-3,\ldots,7 & 11 & \text{multiple of }A\\
y-x=-2 & (1+\xi)^{10} & x=-3,\ldots,6 & 10 & 0\\
x=-3 & (1+y)^9 & y=-4,\ldots,4 & 9 & 0\\
y=4 & (1+x)^8 & x=-2,\ldots,5 & 8 & 0\\
y=3 & (1+x)^6 & x=-2,\ldots,4 & 7 & \text{multiple of }B\\
x=-2 & (1+y)^6 & y=-3,\ldots,2 & 6 & 0\\
y-x=-1 & (1+\xi)^5 & x=-1,\ldots,3 & 5 & 0\\
y=2 & (1+x)^4 & x=-1,\ldots,2 & 4 & 0\\
x=-1 & (1+y)^3 & y=-1,0,1 & 3 & 0\\
y=1 & (1+x)^2 & x=0,1 & 2 & 0
\end{array}
\]
The table is processed from top to bottom.  Its first row determines a scalar multiple of the top slice of $A$.  Subtract that multiple of $A$, which belongs to $V$ by Lemma~\ref{lem:edge-mutability}.  The next three rows then vanish in order by the stated divisibility and exponent-count criterion.  The fifth row determines a scalar multiple of the $y=3$ slice of $B$, namely $x^{-2}(1+x)^6$.  Subtract that multiple of $B$, which also belongs to $V$.  Since the support of $B$ is contained in
\[
\operatorname{Supp}(B)\subseteq\{x\ge-2,\ y\le3,\ y-x\ge-1\},
\]
this subtraction leaves the first four eliminated slices unchanged.  The remaining rows then vanish successively, since each has at most as many consecutive exponents as its divisibility exponent.

After these eliminations the remaining support satisfies $y\le0$, $x\ge0$, and $y-x\ge0$, so it is supported only at $(0,0)$.  Thus every element of $V$ belongs to $\operatorname{span}\{A,B,1\}$.  Lemma~\ref{lem:edge-mutability} shows that each of $A$, $B$, and $1$ belongs to $V$, proving the reverse inclusion.
\end{proof}

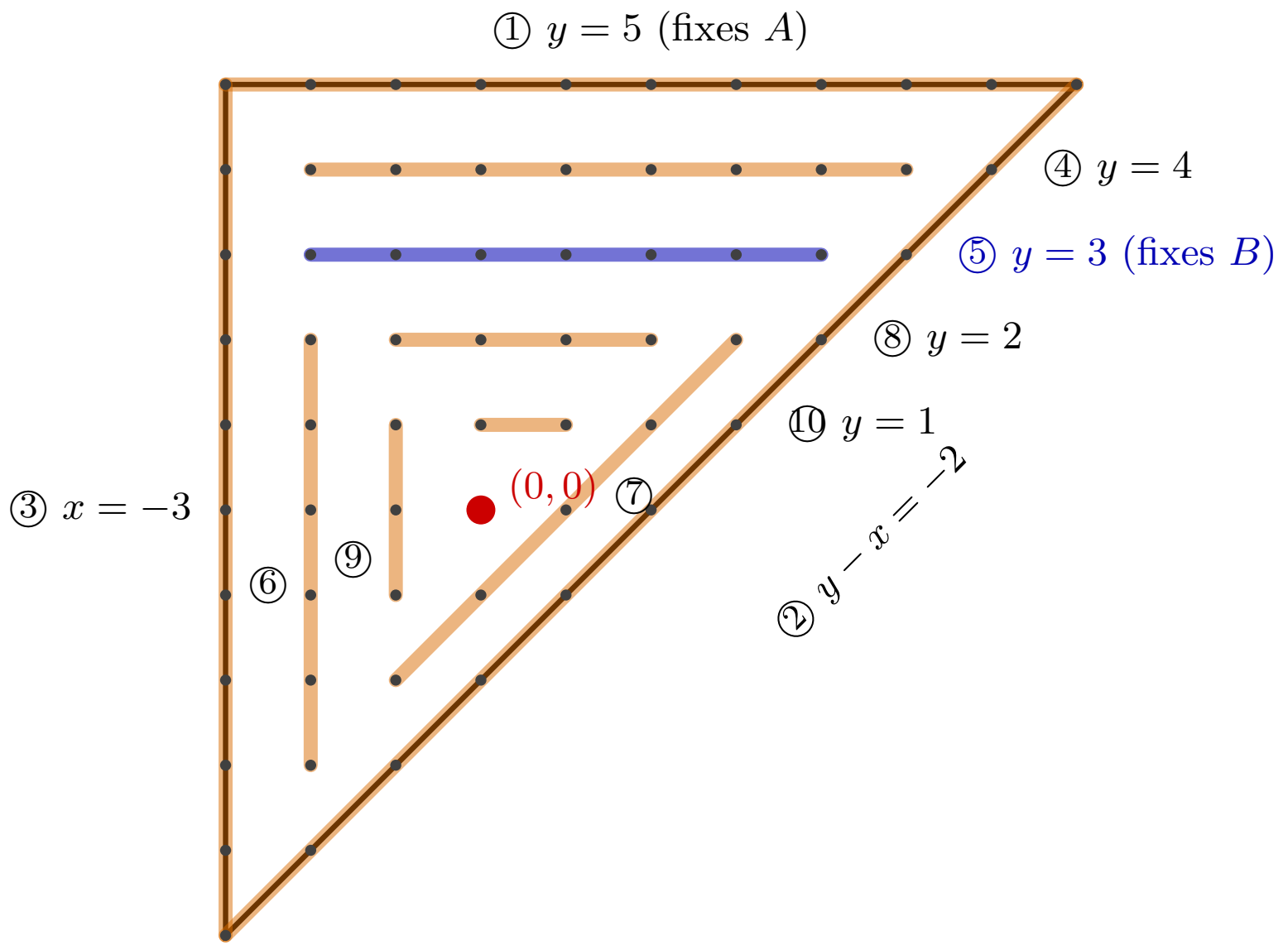
\begin{figure}[!htbp]
\centering
\begin{tikzpicture}[scale=0.68]
  \draw[very thick] (-3,-5) -- (-3,5) -- (7,5) -- cycle;
  \tikzset{sl/.style={orange!85!black,line width=3.2pt,opacity=.5,line cap=round}}
  \draw[sl] (-3,5) -- (7,5);            % 1
  \draw[sl] (-3,-5) -- (7,5);           % 2
  \draw[sl] (-3,-5) -- (-3,5);          % 3
  \draw[sl] (-2,4) -- (5,4);            % 4
  \draw[blue!70!black,line width=3.2pt,opacity=.55,line cap=round] (-2,3) -- (4,3);  % 5
  \draw[sl] (-2,-3) -- (-2,2);          % 6
  \draw[sl] (-1,-2) -- (3,2);           % 7
  \draw[sl] (-1,2) -- (2,2);            % 8
  \draw[sl] (-1,-1) -- (-1,1);          % 9
  \draw[sl] (0,1) -- (1,1);             % 10
  \foreach \x in {-3,...,7}
    \foreach \y in {-5,...,5}{
      \pgfmathtruncatemacro{\d}{\y-\x}
      \ifnum\d>-3 \fill[black!75] (\x,\y) circle (0.065); \fi
    }
  \fill[red!80!black] (0,0) circle (0.17);
  \node[red!80!black,font=\footnotesize,right] at (0.12,0.25) {$(0,0)$};
  \node[font=\footnotesize,above] at (2,5.2) {\textcircled{\scriptsize 1} $y=5$ (fixes $A$)};
  \node[font=\footnotesize,rotate=45] at (4.6,-0.4) {\textcircled{\scriptsize 2} $y-x=-2$};
  \node[font=\footnotesize,left] at (-3.2,0) {\textcircled{\scriptsize 3} $x=-3$};
  \node[font=\footnotesize,right] at (6.4,4) {\textcircled{\scriptsize 4} $y=4$};
  \node[font=\footnotesize,right,blue!70!black] at (5.4,3) {\textcircled{\scriptsize 5} $y=3$ (fixes $B$)};
  \node[font=\footnotesize] at (-2.5,-0.9) {\textcircled{\scriptsize 6}};
  \node[font=\footnotesize] at (1.8,0.15) {\textcircled{\scriptsize 7}};
  \node[font=\footnotesize,right] at (4.4,2) {\textcircled{\scriptsize 8} $y=2$};
  \node[font=\footnotesize] at (-1.5,-0.6) {\textcircled{\scriptsize 9}};
  \node[font=\footnotesize,right] at (3.4,1) {\textcircled{\scriptsize 10} $y=1$};
\end{tikzpicture}
\caption{The elimination in the proof of Lemma~\ref{lem:direct-T-space}. Slices are removed in the numbered order (\textcircled{\scriptsize 6}: $x=-2$, \textcircled{\scriptsize 7}: $y-x=-1$, \textcircled{\scriptsize 9}: $x=-1$). Steps \textcircled{1} and \textcircled{5} each determine one scalar, a multiple of $A$ and of $B$. Every other slice is forced to vanish by the exponent count, and after the ten steps only the origin remains, so the surviving space is $\operatorname{span}\{A,B,1\}$.}
\label{fig:elimination}
  \end{figure}
\FloatBarrier

\begin{proof}[Proof of Theorem~\ref{thm:mmlp-family}]
By Lemma~\ref{lem:direct-T-space}, an element of $V$ has the form $c_AA+c_BB+c_0$.  Vertex normalization forces $c_A=1$, and the zero-constant-term condition gives $c_0=-2520-60c_B$ by \eqref{eq:basic-CT}.  Since $B$ is supported strictly inside $P$, these polynomials have the fixed edge profiles of $A$, which are the prescribed $T$-binomial coefficients for this polygon.  Writing $s=c_B$, we obtain $f_s=A+sB-(2520+60s)$.  Since $[x^4y^3]A=360$ and $[x^4y^3]B=1$ by \eqref{eq:basic-a}, its distinguished coefficient is $a=s+360$.  This gives exactly $f_a$ in \eqref{eq:fa}.  Proposition~\ref{prop:mmlp-facts}(1) identifies these polynomials with the MMLPs.
\end{proof}

The distinguished coefficient also records the residual point.  On the left edge, the decomposition $10=3\cdot3+1$ leaves the residual cone
\[
R=\operatorname{conv}\{0,(-3,4),(-3,5)\}.
\]
The point $(-2,3)$ lies in its interior and is the residual point in the CKPT construction \cite[Definitions 3.5--3.6]{CKPT}.  Direct expansion gives
\[
[x^{-2}y^3]A=[x^4y^3]A=360,\qquad [x^{-2}y^3]B=[x^4y^3]B=1.
\]
Thus $a_{[4,3]}=[x^4y^3]f_a=[x^{-2}y^3]f_a$, so Oneto--Petracci's distinguished coordinate is the coefficient at this residual point in our family \cite[Section 2]{OP}.

\begin{lemma}\label{lem:parameter-dimensions}
For polygon n.1 and $\cX_{10}$, both parameter spaces are one-dimensional:
\[
\dim L_P^T=\dim H^{\mathrm{ts}}_{\cX_{10}}=1.
\]
\end{lemma}

\begin{proof}
The twisted side is one-dimensional because $\alpha$ is the unique sector of age less than one.  Lemma~\ref{lem:direct-T-space} gives $V=\operatorname{span}\{A,B,1\}$.  Vertex normalization fixes the coefficient of $A$, and the zero-constant-term condition fixes the constant coefficient, leaving the coefficient of $B$ free.  Hence $\dim L_P^T=1$.
\end{proof}

With the complete family and its one-dimensionality established, we now turn to a general fiber $f_s(x,y)=\eta$.  Set $T=\eta+C(s)$, so its equation is
\begin{equation}\label{eq:fiberAB}
A+sB=T.
\end{equation}

Since $A$ and $B$ share the factor $H^6$, the following substitutions simplify the fiber equation.
\begin{equation}\label{eq:uz-def}
z=\frac{yH^2}{x},
\qquad
u=xz^2=\frac{y^2H^4}{x}.
\end{equation}
Indeed, the definition of $z$ gives $B=xz^3$.  The further substitution $u=xz^2$ simplifies both summands simultaneously:
\[
B=uz,
\qquad
A=u^2z.
\]
Thus the fiber equation becomes linear in $z$:
\[
zu(u+s)=T.
\]
Solving for $z$ and using the relation $u=xz^2$ gives
\begin{equation}\label{eq:z-x-nu}
z=\frac{T}{u(u+s)},
\qquad
x=\frac{u}{z^2}=\frac{u^3(u+s)^2}{T^2}.
\end{equation}

The defining relation for $z$, meanwhile, is quadratic in $y$:
\[
z=\frac{((1+x)y+1)^2}{xy}.
\]
Therefore $y$ satisfies the quadratic equation
\begin{equation}\label{eq:y-quadratic}
(1+x)^2y^2+\bigl(2(1+x)-xz\bigr)y+1=0.
\end{equation}
The quadratic polynomial has discriminant
\[
\Delta=xz\bigl(xz-4(1+x)\bigr).
\]
The auxiliary variable $w=2(1+x)^2y+2(1+x)-xz$ therefore satisfies $w^2=\Delta$.  Substitution of \eqref{eq:z-x-nu} yields
\[
w^2=
\frac{u^2(u+s)}{T^3}
\left(Tu^2(u+s)-4T^2-4u^3(u+s)^2\right).
\]
After the rescaling
\[
v=\frac{\sqrt{-1}\,T^2}{u}w,
\]
we obtain \eqref{eq:hyperelliptic}.

\begin{proof}[Proof of Theorem~\ref{thm:hyperelliptic}]
On the open set where the product $Tu(u+s)(1+x)$ is nonzero, a point $(u,v)$ of \eqref{eq:hyperelliptic} determines the inverse coordinates
\begin{equation}\label{eq:inverse-birational}
z=\frac{T}{u(u+s)},
\qquad
x=\frac{u^3(u+s)^2}{T^2},
\qquad
w=-\frac{\sqrt{-1}\,u}{T^2}v.
\end{equation}
The quadratic equation recovers $y$ as
\[
y=\frac{w-2(1+x)+xz}{2(1+x)^2}.
\]
Together with \eqref{eq:uz-def}, these formulas are inverse on their domains.  When $T$ is nonzero, the fiber relation $zu(u+s)=T$ excludes $u=0$ and $u=-s$, and these values account for only finitely many points on the hyperelliptic curve.  The condition $x=-1$ in the inverse formulas is equivalent to
\[
u^3(u+s)^2+T^2=0,
\]
which has only finitely many roots.  Each root gives at most two points on \eqref{eq:hyperelliptic}, and the defining equation for $y$ is linear with nonzero coefficient $z$.  Thus only finitely many points are omitted, and the rational maps identify nonempty open subsets of the two curves, proving birationality.

The right-hand side of \eqref{eq:hyperelliptic} has leading term $4T u^6$, so it has degree exactly six when $T$ is nonzero.  Squarefreeness is open in $(s,T)$, so it suffices to exhibit one squarefree specialization.  At $s=0$ and $T=1$, the sextic becomes
\[
u\bigl(4u^5-u^3+4\bigr).
\]
The factor $4u^5-u^3+4$ does not vanish at $u=0$.  A repeated nonzero root would satisfy
\[
4u^5-u^3+4=0,
\qquad
20u^4-3u^2=0.
\]
The second equation forces $u^2=3/20$.  Substitution in the first yields $u^3=10$, from which it follows that $u^6$ equals $100$.  But $u^2=3/20$ also forces $u^6=(3/20)^3$, a contradiction.  This proves that the sextic is squarefree at this specialization, so its discriminant is not identically zero.  For general $(s,T)$ it is squarefree and hence not a square in $\C(u)$, so the hyperelliptic curve is irreducible.  The open-set correspondence above then proves that the original fiber is irreducible.  A smooth double cover of $\Pj^1$ branched at six points has genus two.
\end{proof}

\section{An \texorpdfstring{$S$}{S}-extended GIT presentation of \texorpdfstring{$X_{10}$}{X10}}\label{sec:GIT}

To compute the quantum period of $\cX_{10}$, we construct an $S$-extended GIT presentation and verify that its data satisfy the hypotheses of Wang's mirror theorem, which does not require convexity.  The usual quantum Lefschetz argument would relate its genus-zero Gromov--Witten theory to that of the ambient stack, but it requires convexity of the defining line bundle \cite[Introduction]{Wang}: $H^1(C,f^*L)=0$ for every genus-zero stable map $f$.  Convexity is not automatic here.  Since $10$ is not divisible by the weight $3$, $\mathcal O(10)$ has nontrivial isotropy along the $\mu_3$-locus, so we use the explicit $I$-function of the extended presentation for the direct calculation.

The ordinary weighted-projective presentation does not attach an independent formal variable to the unique age-$\frac23$ sector in Conjecture~B.  We therefore use the $S$-extension associated with its box element, adding an auxiliary coordinate $u$ and a second torus factor.  The new effective-degree direction records the twisted-sector deformation in the $I$-function.  Wang's Section~7 example uses an extended GIT construction for a cubic hypersurface in $\Pj(1,1,1,2)$ to compute its small quantum product with twisted-sector data \cite[Section 7]{Wang}, providing a precedent for the present construction.

Extended GIT presentations, the GIT analogue of $S$-extensions, are established tools for nonconvex complete intersections.  Janda, Sultani, and Zhou develop them for complete intersections in weighted projective stacks and show that their extended $I$-function agrees with Wang's big $I$-function for an appropriate GIT extension \cite{JandaSultaniZhou}.  We now fix the GIT and quasimap conventions needed to apply Wang's theorem to the $X_{10}$ presentation.

Let a torus $G$ act linearly on
\[
W=\bigoplus_{\rho=1}^n\C_\rho,
\]
and let $\theta$ be a stability character.  When stable and semistable loci agree, the quotient is
\[
\mathcal T=[W^{ss}(\theta)/G].
\]
A character $\xi$ of $G$ determines a line bundle $L_\xi$ on $\mathcal T$.  A quasimap degree is a homomorphism
\[
\beta\in\operatorname{Hom}_{\Z}\bigl(\operatorname{Pic}([W/G]),\Q\bigr),
\]
which records the rational numbers $\beta(L_\xi)$, and Wang's effective semigroup is denoted by
\[
\Eff(W,G,\theta)
\]
\cite[Definition 2.4]{Wang}.  It consists of finite sums of degree classes represented by quasimaps to the quotient.  In a character basis $\pi_1,\ldots,\pi_k$ of $G$, the fractional parts determine the torsion element
\[
g_\beta
=\bigl(e^{2\pi\sqrt{-1}\,\beta(L_{\pi_1})},\ldots,
e^{2\pi\sqrt{-1}\,\beta(L_{\pi_k})}\bigr)\in G,
\]
which determines the fixed locus used in the construction.  In Wang's $I$-function, however, the $\beta$-summand lies in the inertia component indexed by $g_\beta^{-1}$, not $g_\beta$ \cite[Definition 3.5 and Section 3.1, especially (3.8)]{Wang}.  Wang's distinguished fixed locus is
\begin{equation}\label{eq:Wang-fixed-locus}
Z_\beta^{ss}
=(W^{ss}(\theta))^{g_\beta}\cap
\{x_\rho=0:\beta(L_\rho)\in\Z_{<0}\}.
\end{equation}
In particular, if $\beta(L_\rho)$ is a negative integer, then the coordinate $x_\rho$ is forced to vanish on $Z_\beta^{ss}$ \cite[Section 3, (3.1), and Definition 3.3]{Wang}.

\begin{theorem}[Wang \wangmirrorcite, specialized restatement]\label{thm:wang-specialized}
Let $W=\A^n$ carry a linear action of a torus $G=(\C^\times)^k$, and let $\theta$ be a character such that the stable and semistable loci coincide.  Define the quotient stack $\mathcal T$ by
\[
\mathcal T=[W^{ss}(\theta)/G]
\]
and suppose that $\mathcal T$ is a smooth proper toric Deligne--Mumford stack.  Let $E=\bigoplus_{b=1}^cL_{\tau_b}$ be a direct sum of line bundles associated with characters of $G$, and suppose that $G$-equivariant sections of these bundles cut out an irreducible complete intersection $Z\subset W$ that is smooth along $W^{ss}(\theta)$.  Define the quotient stack $\mathcal Y$ by
\[
\mathcal Y=[Z\cap W^{ss}(\theta)/G].
\]
The series
\[
-zI_{\mathcal Y}(q,\mathbf t,-z)
\]
is a slice of Givental's Lagrangian cone of $\mathcal Y$, where $I_{\mathcal Y}$ is Wang's big $I$-function from \cite[Definition 3.5]{Wang}.  More explicitly, the mirror transformation
\[
\mu(q,\mathbf t,z)
=\bigl[zI_{\mathcal Y}(q,\mathbf t,z)-z\1_{\mathcal Y}\bigr]_+
\]
satisfies
\[
I_{\mathcal Y}(q,\mathbf t,z)
=J_{\mathcal Y}(q,\mu(q,\mathbf t,z),z),
\]
where $[\ ]_+$ denotes truncation to nonnegative powers of $z$.
\end{theorem}

Throughout this paper, the \emph{normalized $I$-function} means the specialization $I_{\mathcal Y}(q,0,z)$ with its degree-zero term fixed to $\1_{\mathcal Y}$.  In the present application its expansion begins
\[
I=\1+O(z^{-1}).
\]
Wang's theorem extends the toric-stack mirror theorem of Coates--Corti--Iritani--Tseng \cite{CCIT}, and the weighted-projective $J$-function conventions used below are compatible with \cite{CCLT}.

For the age-$\frac23$ sector, set $\zeta:=e^{2\pi\sqrt{-1}/3}$.  Its box element is $g=\zeta^2=\zeta^{-1}$, whose fractional parts on the ambient weighted-projective coordinates are $(\frac23,\frac13,0,\frac13)$.  The $S$-extension assigns the auxiliary weights $\lfloor\frac23(1,2,3,5)\rfloor=(0,1,2,3)$ to $x_0,x_1,x_2,x_3$ and weight $1$ to $u$.
Accordingly, the torus $G=(\C^\times)^2$ acts on $W=\A^5$ with charge matrix
\begin{equation}\label{eq:charge}
Q^S=
\begin{pmatrix}
1&2&3&5&0\\
0&1&2&3&1
\end{pmatrix}.
\end{equation}

The two rows of \eqref{eq:charge} determine the standard characters $\chi_t=(1,0)$ and $\chi_u=(0,1)$.  Let $L_t$ and $L_u$ be the associated line bundles on the quotient, and for a degree class $\beta$ write $d_t=\beta(L_t)$ and $d_u=\beta(L_u)$.

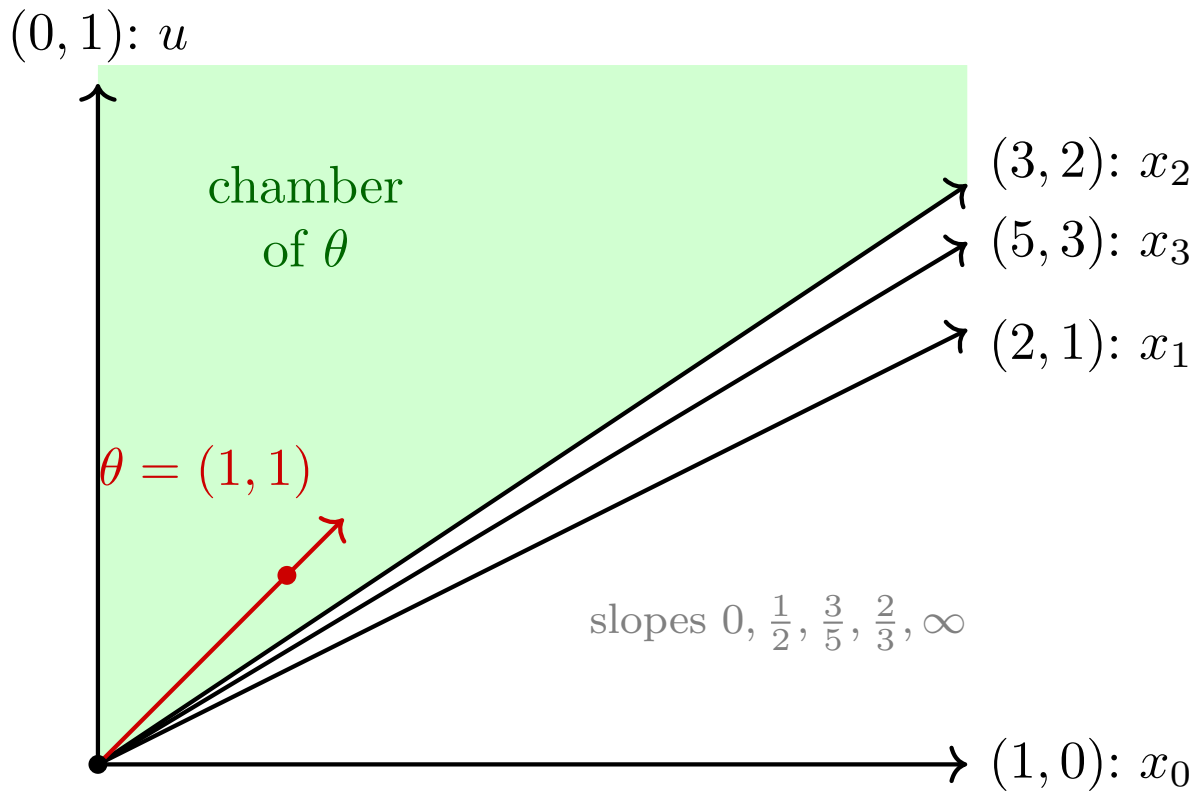
\begin{figure}[htbp]
\centering
\begin{tikzpicture}[scale=1.25]
  \fill[green!18] (0,0) -- (4.6,3.067) -- (4.6,3.7) -- (0,3.7) -- cycle;
  \draw[thick,->] (0,0) -- (4.6,0)      node[right,font=\small] {$(1,0)$:\ $x_0$};
  \draw[thick,->] (0,0) -- (4.6,2.3)    node[right,font=\small,yshift=-3pt] {$(2,1)$:\ $x_1$};
  \draw[thick,->] (0,0) -- (4.6,2.76)   node[right,font=\small] {$(5,3)$:\ $x_3$};
  \draw[thick,->] (0,0) -- (4.6,3.067)  node[right,font=\small,yshift=4pt] {$(3,2)$:\ $x_2$};
  \draw[thick,->] (0,0) -- (0,3.6)      node[above,font=\small] {$(0,1)$:\ $u$};
  \draw[red!80!black,thick,->] (0,0) -- (1.3,1.3);
  \fill[red!80!black] (1,1) circle (0.05);
  \node[red!80!black,font=\small,above left] at (1.25,1.3) {$\theta=(1,1)$};
  \node[green!40!black,font=\small,align=center] at (1.1,2.9) {chamber\\ of $\theta$};
  \node[font=\scriptsize,gray] at (3.6,0.75) {slopes $0,\tfrac12,\tfrac35,\tfrac23,\infty$};
  \fill (0,0) circle (0.05);
\end{tikzpicture}
\caption{The weights of the $S$-extended presentation \eqref{eq:charge} (columns of $Q^S$) and the stability character $\theta=(1,1)$, which lies in the open chamber bounded by the rays through $(3,2)$ and $(0,1)$. A point is semistable exactly when its nonzero coordinates span a cone containing $\theta$ in its interior, which forces $u\ne0$ and some $x_i\ne0$.}
\label{fig:chamber}
\end{figure}

We choose the character $\theta$ to have coordinates $(1,1)$.  The columns of \eqref{eq:charge} are
\[
(1,0),(2,1),(3,2),(5,3),(0,1).
\]
In increasing order of slope, which differs from the displayed column order, their rays have slopes
\[
0,\quad \frac12,\quad \frac35,\quad \frac23,\quad \infty.
\]
Since the chosen character has slope $1$, it lies in the open chamber between the rays generated by $(3,2)$ and $(0,1)$.  By the torus cone criterion, semistability requires both $u\ne0$ and at least one $x_i\ne0$: without $u$, all available weights have slope at most $\frac23$, while without the $x_i$ only the ray $(0,1)$ remains.  Conversely, the weights of $u$ and any nonzero $x_i$ span a cone containing $\theta$ in its interior and have finite stabilizer, so every such point is stable.  Hence the stable and semistable loci coincide and are described by
\[
W^{ss}(\theta)=W^s(\theta)=
\{u\ne0,\ (x_0,x_1,x_2,x_3)\ne0\}.
\]
We use the second $\C^\times$ to normalize $u=1$, after which the residual action has weights $1,2,3,5$.  Therefore the quotient stack is isomorphic to the weighted projective stack
\begin{equation}\label{eq:quotient}
[W^{ss}(\theta)/G]\simeq\Pj(1,2,3,5).
\end{equation}
Thus the extension enlarges the degree bookkeeping without changing the ambient toric Deligne--Mumford stack.  Figure~\ref{fig:chamber} displays the ordered weight rays and the chamber of $\theta$.

\FloatBarrier

Wang's theorem treats the hypersurface as a GIT complete intersection, so its defining equation must be equivariant for both torus factors.  We write a weighted-homogeneous polynomial of degree $10$ as
\[
F_{10}=\sum_M c_M M
\]
and express each monomial as $M=x_0^{a_0}x_1^{a_1}x_2^{a_2}x_3^{a_3}$.  The second grading assigns the degree $\nu(M)=a_1+2a_2+3a_3$ to this monomial.  For each weighted-degree-$10$ monomial, $3\nu(M)\le2\deg_{\mathrm{wt}}(M)=20$, since the variable weights $(1,2,3,5)$ and second degrees $(0,1,2,3)$ satisfy this inequality termwise.  As $\nu(M)$ is integral, $\nu(M)\le6$.  The homogenized polynomial
\begin{equation}\label{eq:extended-F}
\widetilde F=\sum_M c_M M u^{6-\nu(M)}
\end{equation}
has nonnegative exponents of $u$, and every term has the common $G$-character $(10,6)$.  On the stable locus, normalizing $u=1$ recovers $F_{10}$.

\begin{lemma}[Applicability of Wang's mirror theorem]\label{lem:Wang-hypotheses}
For a general weighted-homogeneous polynomial $F_{10}$ of degree $10$, the $S$-extended data $(W,G,\theta,\widetilde F)$ satisfy the hypotheses of the specialized restatement in \cref{thm:wang-specialized}.  Explicitly, $W$ is the smooth irreducible affine l.c.i.\ representation $\A^5$ of $G=(\C^\times)^2$, the equality $W^{ss}(\theta)=W^s(\theta)$ holds, and the quotient is the smooth proper toric Deligne--Mumford stack $\Pj(1,2,3,5)$.  Moreover, $\widetilde F$ is $G$-equivariant of character $(10,6)$, and its zero locus is irreducible and smooth along $W^s(\theta)$.
\end{lemma}

\begin{proof}
The affine space $W=\A^5$ is a smooth irreducible l.c.i.\ representation of the torus $G=(\C^\times)^2$.  The chamber calculation gives $W^s(\theta)=W^{ss}(\theta)$, and \eqref{eq:quotient} identifies the quotient with the smooth proper toric stack $\Pj(1,2,3,5)$.

Every term of $\widetilde F$ has $G$-degree $(10,6)$, so $\widetilde F$ is a $G$-equivariant section of the line bundle associated with that character.  On $u\ne0$, gauge-fixing $u=1$ identifies its zero locus with the affine cone of $F_{10}$ away from the origin.  For example, the weighted-homogeneous polynomial
\[
F_{10}^{\circ}=x_0^{10}+x_1^5+x_0x_2^3+x_3^2
\]
has partial derivatives $10x_0^9+x_2^3$, $5x_1^4$, $3x_0x_2^2$, and $2x_3$.  Their common vanishing forces $x_1=x_3=0$ and $x_0x_2^2=0$, after which $10x_0^9+x_2^3=0$ forces $x_0=x_2=0$.  Thus this degree-$10$ hypersurface is quasismooth.  Since quasismoothness is open in the space of weighted-homogeneous degree-$10$ polynomials, a general $F_{10}$ is quasismooth and gives smoothness over $W^s(\theta)$.

Finally, quasismoothness forces $F_{10}$ to be irreducible.  If it had two distinct positive-degree factors, their hypersurfaces in $\A^4$ would intersect in dimension at least two, so they would meet away from the origin.  At an intersection point the gradient of $F_{10}$ vanishes, contradicting quasismoothness.  If $F_{10}$ had a repeated factor, its gradient would vanish along that factor's hypersurface, which contains nonzero points, again contradicting quasismoothness.  Since $G$ is a torus, it preserves each irreducible factor of the semi-invariant polynomial $\widetilde F$ up to scalar, so each factor is itself semi-invariant.  If $\widetilde F$ factored into nonconstant factors, setting $u=1$ would factor $F_{10}$ unless one factor became constant.  Such a factor has first $G$-degree zero, so it is a polynomial in $u$ alone.  Its semi-invariance under the second torus makes it a monomial in $u$.  But $u$ does not divide $\widetilde F$: the general polynomial contains $x_3^2$ with nonzero coefficient, of extended degree $2(5,3)=(10,6)$, so this term occurs with $u^0$.  Hence $\widetilde F$ is irreducible.
\end{proof}

\section{The quantum period and proof of Conjecture B}\label{sec:Ifunction}\label{sec:period-compare}

We compute the quantum period of $\cX_{10}$ from the $S$-extended $I$-function and prove Theorem~\ref{thm:conjB}.  We begin by determining the effective degrees whose inertia components meet $\cX_{10}$.  The only stacky coordinate points of $\Pj(1,2,3,5)$ have stabilizers $\mu_2$, $\mu_3$, and $\mu_5$.  A general degree-$10$ polynomial $F_{10}$ contains the terms $x_1^5$ and $x_3^2$ with nonzero coefficients, so its zero locus misses the $\mu_2$- and $\mu_5$-coordinate points.  By contrast, the hypersurface contains the $\mu_3$-coordinate point because no pure power of $x_2$ has weighted degree $10$.

\begin{lemma}[Relevant effective degrees]\label{lem:effective-semigroup}
The classes in $\Eff(W,G,\theta)$ whose associated inertia components, indexed by $g_\beta^{-1}$, meet $\cX_{10}$ are exactly
\[
\mathbb N\beta_1+\mathbb N\beta_2.
\]
The two generators are characterized by
\begin{equation}\label{eq:betas}
\bigl(\beta_1(L_t),\beta_1(L_u)\bigr)=\left(\frac13,0\right),
\qquad
\bigl(\beta_2(L_t),\beta_2(L_u)\bigr)=\left(-\frac23,1\right).
\end{equation}
\end{lemma}

\begin{proof}
By Wang's definition, an effective class is a finite sum of degrees represented by quasimaps \cite[Definition 2.4]{Wang}.  For each quasimap degree, write $d_t,d_u$ for its pairings with $L_t,L_u$.  A quasimap must map the generic point of its source curve into the semistable locus.  In this chamber, that requires $u$ and at least one of $x_0,x_1,x_2,x_3$ to be generically nonzero.  The corresponding sections have nonnegative degrees, so $d_u\ge0$ and at least one of the coordinate degrees
\[
d_0=d_t,\qquad d_1=2d_t+d_u,\qquad d_2=3d_t+2d_u,\qquad d_3=5d_t+3d_u
\]
is nonnegative.  The identities
\[
3d_t+2d_u=3d_0+2d_u=\frac32d_1+\frac12d_u=d_2=\frac35d_3+\frac15d_u
\]
then show that $3d_t+2d_u\ge0$ for each quasimap degree.  By additivity, both inequalities hold for every effective class.

Only the identity and $\mu_3$ inertia sectors meet $\cX_{10}$.  For a relevant effective class $\beta$, write $d_t=\beta(L_t)$ and $d_u=\beta(L_u)$.  In these sectors the inertia elements act trivially on $u$ and $x_2$.  By the definition of $g_\beta$, the corresponding degree pairings
\[
k:=d_u,\qquad m:=\beta(L_{x_2})=3d_t+2d_u
\]
are integers.  The inequalities give $m,k\ge0$, so the degree pairings have the form
\[
(d_t,d_u)=m\left(\frac13,0\right)+k\left(-\frac23,1\right).
\]
This expression is unique.  Here $k$ is the degree of the extension coordinate, whereas $m$ is the $x_2$-degree.

Conversely, representable quasimaps realize both generators using Wang's description of toric quasimaps by character line bundles and coordinate sections \cite[Example 2.6]{Wang}.  We take
\[
\pi:\mathcal C=\sqrt[3]{(\Pj^1,p)}\longrightarrow\Pj^1
\]
to be the third root stack at a point $p$, and we denote its tautological root line bundle by $R$.  The root line bundle satisfies
\[
R^3\simeq\pi^*\mathcal O_{\Pj^1}(p),
\qquad
\deg R=\frac13.
\]
For brevity, we write $\mathcal O_{\mathcal C}(p)=\pi^*\mathcal O_{\Pj^1}(p)$.
Write $\mathcal L_t$ and $\mathcal L_u$ for the pullbacks to $\mathcal C$ of the quotient line bundles $L_t$ and $L_u$.

For $\beta_1$, we keep the extension coordinate nonzero and let $x_2$ carry the base point.  We take
\[
\mathcal L_t=R,
\qquad
\mathcal L_u=\mathcal O_{\mathcal C}.
\]
The $x_2$-coordinate is therefore a section of
\[
\mathcal L_t^3\mathcal L_u^2\simeq\mathcal O_{\mathcal C}(p).
\]
We choose $u$ to be a nonzero constant and $x_2$ to be a section with a single zero at an ordinary point away from $p$, setting all other coordinates to zero.  Away from that zero, both $x_2$ and $u$ are nonzero, so the image lies in $W^{ss}(\theta)$.  The zero of $x_2$ is the unique, isolated base point.  At the order-three stacky point, the stabilizer acts on $R$ through the standard character, so the induced homomorphism $\mu_3\to G$ is $\zeta\mapsto(\zeta,1)$ and is injective.  The quasimap is therefore representable and has degree
\[
\bigl(\deg\mathcal L_t,\deg\mathcal L_u\bigr)=\left(\frac13,0\right)=\beta_1.
\]

For $\beta_2$, the roles are reversed: $x_2$ remains nonzero and the extension coordinate carries the base point.  We take
\[
\mathcal L_t=R^{-2},
\qquad
\mathcal L_u=\mathcal O_{\mathcal C}(p).
\]
These choices give
\[
\mathcal L_t^3\mathcal L_u^2
\simeq\mathcal O_{\mathcal C}(-2p)\otimes\mathcal O_{\mathcal C}(2p)
\simeq\mathcal O_{\mathcal C},
\]
so we choose $x_2$ to be a nonzero constant.  We choose $u$ to be a section of $\mathcal O_{\mathcal C}(p)$ with a single zero at an ordinary point away from $p$, setting all other coordinates to zero.  The generic image lies in $W^{ss}(\theta)$, and the zero of $u$ is the unique, isolated base point.  The stabilizer character of $R^{-2}$ is $\zeta^{-2}=\zeta$, so the induced homomorphism at the stacky point is again $\zeta\mapsto(\zeta,1)$.  Hence this quasimap is representable and has degree
\[
\bigl(\deg\mathcal L_t,\deg\mathcal L_u\bigr)=\left(-\frac23,1\right)=\beta_2.
\]

Thus both $\beta_1$ and $\beta_2$ belong to $\Eff(W,G,\theta)$.  Since $\Eff(W,G,\theta)$ is a semigroup, every $m\beta_1+k\beta_2$ with $m,k\ge0$ is effective.  Together with the necessary form above, this proves the statement about the effective semigroup.  It does not assert that $Z_\beta^{ss}$ is nonempty or that the localized Gysin contribution is nonzero.
\end{proof}

For an effective degree $\beta=m\beta_1+k\beta_2$ with $m,k\ge0$, the five coordinate degrees are
\begin{equation}\label{eq:degrees}
\begin{aligned}
d_0&=\frac{m-2k}{3},&
 d_1&=\frac{2m-k}{3},\\
d_2&=m,&
 d_3&=\frac{5m-k}{3},\\
d_u&=k.
\end{aligned}
\end{equation}
The hypersurface degree is
\begin{equation}\label{eq:Fdegree}
d_F=\frac{10m-2k}{3}.
\end{equation}

For $d\in\Q$ and a cohomology class $D$, define
\[
\begin{aligned}
R(d,D,z)&=\begin{cases}
\displaystyle\left(\prod_{\substack{i\in\Z\\0\le i<d}}(D+(d-i)z)\right)^{-1},&d>0,\\
    \displaystyle\prod_{\substack{i\in\Z\\d<i<0}}(D+(d-i)z),&d<0,\\
1,&d=0,
\end{cases}\\
R_F(d,D,z)&=\begin{cases}
\displaystyle\prod_{\substack{i\in\Z\\0\le i<d}}(D+(d-i)z),&d>0,\\
\displaystyle\left(\prod_{\substack{i\in\Z\\d<i<0}}(D+(d-i)z)\right)^{-1},&d<0,\\
1,&d=0.
\end{cases}
\end{aligned}
\]

Here $R$ is the coordinate factor and $R_F$ the hypersurface factor in Wang's formula.  Put $D_\rho=c_1(L_\rho)$ for $\rho\in\{0,1,2,3,u\}$ and $D_F=c_1(L_{(10,6)})$.  For $\beta=m\beta_1+k\beta_2$, define $q_1=q^{\beta_1}$ and $q_2=q^{\beta_2}$, so its Novikov monomial is $q^\beta=q_1^mq_2^k$, and set
\[
\mathcal G_{m,k}:=i_{\beta*}s^!_{E_\beta,\mathrm{loc}}
\left(\left[Z_\beta^{ss}/(G/\langle g_\beta^{-1}\rangle)\right]\right)
\in H^*(\bar I_{g_\beta^{-1}}\cX_{10}),
\]
where $E_\beta=L_{(10,6)}$ if $d_F\in\Z_{\ge0}$ and is the rank-zero bundle otherwise, and $i_\beta$ is the natural inclusion of the corresponding hypersurface fixed locus into the rigidified inertia component indexed by $g_\beta^{-1}$.  This is Wang's localized Gysin class, including when $d_F<0$.

With the degrees in \eqref{eq:degrees}--\eqref{eq:Fdegree}, the specialization of Wang's normalized $I$-function at zero input is
\begin{equation}\label{eq:specialized-I}
I(q_1,q_2,z)=\sum_{m,k\ge0}q_1^mq_2^k\,R_F(d_F,D_F,z)
\left(\prod_{\rho\in\{0,1,2,3,u\}}R(d_\rho,D_\rho,z)\right)\mathcal G_{m,k}.
\end{equation}

All products in $R$ and $R_F$ are over integers in the indicated ranges.  In particular, the localized Gysin class remains in every summand rather than being replaced by an ambient hypergeometric fundamental class.  This is the hypersurface specialization of Wang's Definition~3.5 and formula~(3.8) for the charge matrix \eqref{eq:charge} and character $(10,6)$ \cite[Definition 3.5 and Section 3.1, especially (3.8)]{Wang}.

The mirror transformation depends on the nonnegative part of $zI-z\1$, so we classify the effective degrees whose nonzero summands in \eqref{eq:specialized-I} begin at order $z^{-1}$.  The leading exponents follow from the displayed product ranges.
\begin{lemma}[$z^{-1}$-degree classification]\label{lem:z-order}
Let $(m,k)\ne(0,0)$ index an effective degree whose localized contribution is nonzero, and let $\nu(m,k)$ be the leading $z$-exponent of its $I$-function summand.  Then $\nu(m,k)=-1$ occurs exactly for
\[
(m,k)\in\{(1,0),(0,1),(3,0),(2,1)\}.
\]
Every other such summand begins at order $z^{-2}$ or lower.
\end{lemma}

\begin{proof}
Taking the leading power of $z$ in the coordinate factor $R(d,D,z)$ of \eqref{eq:specialized-I} gives
\[
c(d)=
\begin{cases}
-\lceil d\rceil,&d\notin\Z_{<0},\\
-d-1,&d\in\Z_{<0},
\end{cases}
\]
with $c(0)=0$.  The leading exponent of the hypersurface factor $R_F(d,D,z)$ is
\[
h(d)=
\begin{cases}
\lceil d\rceil,&d\notin\Z_{<0},\\
d+1,&d\in\Z_{<0}.
\end{cases}
\]
These contributions give the leading exponent
\[
\nu(m,k)=h(d_F)+c(d_0)+c(d_1)+c(d_2)+c(d_3)+c(k),
\]
for the $(m,k)$-summand, where the degrees are given by \eqref{eq:degrees}--\eqref{eq:Fdegree}.  We put $N=m+k$ and write
\[
\nu(m,k)=-\frac{N}{3}+\delta(m,k).
\]
As $(m,k)$ varies, the formula changes only when $d_0$, $d_1$, or $d_3$ crosses zero, giving the walls
\[
2k=m,\qquad k=2m,\qquad k=5m.
\]
The hypersurface degree introduces no further wall because $d_F$ equals $2d_3$.  Within each region $\delta$ depends only on $N$ modulo $3$, as recorded in Table~\ref{tab:wall-regions}.

\begin{table}[htbp]
\centering
\caption{Values of $\delta(m,k)$ in the four wall regions.}\label{tab:wall-regions}
\begin{tabular}{c c c c c}
\toprule
Region & inequalities & $N\bmod3=0$ & $N\bmod3=1$ & $N\bmod3=2$ \\
\midrule
$\mathrm I$ & $2k\le m$ & $0$ & $-\frac23$ & $-\frac43$ \\
$\mathrm {II}$ & $m<2k,\ k\le2m$ & $-1$ & $-\frac23$ & $-\frac43$ \\
$\mathrm {III}$ & $2m<k,\ k\le5m$ & $-2$ & $-\frac23$ & $-\frac43$ \\
$\mathrm {IV}$ & $5m<k$ & $-2$ & $-\frac23$ & $-\frac43$ \\
\bottomrule
\end{tabular}
\end{table}

In Region~I all relevant degrees are nonnegative, giving $\delta=0,-2/3,-4/3$ according to $N\bmod3$.  Crossing $2k=m$ makes $d_0$ negative.  When $N$ is congruent to $0$ modulo $3$, this is a negative integer and lowers $\delta$ by one.  Crossing $k=2m$ has the same effect for $d_1$.  At the final wall, the negative-integral corrections from $d_3$ and $d_F=2d_3$ cancel.  Thus Regions~III and IV have the same row.

For $N\equiv2\pmod3$, Table~\ref{tab:wall-regions} gives $\nu\le-2$.  For $N\equiv1\pmod3$, the equality $\nu=-1$ forces $N=1$, which gives $(m,k)=(1,0)$ or $(0,1)$.  For $N\equiv0\pmod3$, Regions~II--IV give $\nu\le-2$, whereas Region~I gives $\nu=-N/3$.  Consequently, the equality $\nu=-1$ forces $N=3$, and the inequality $2k\le m$ gives $(m,k)=(3,0)$ or $(2,1)$.  These four and only these four degrees contribute at order $z^{-1}$.
\end{proof}

Lemma~\ref{lem:z-order} identifies the four effective degrees contributing at order $z^{-1}$, which are recorded in Table~\ref{tab:zminusone-degrees}.  For degrees $(1,0)$ and $(0,1)$, $g_\beta=(\zeta,1)$ acts on the tangent space of $\cX_{10}$ at the $\frac13(1,1)$ point with weights $(2/3,2/3)$, so $g_\beta$ has age $4/3$.  Its inverse $g_\beta^{-1}=(\zeta^2,1)$ has age $2/3$ and labels the sector $\alpha$ in Wang's $I$-function.  For degrees $(3,0)$ and $(2,1)$, $g_\beta$ is the identity element.

\begin{table}[htbp]
\centering
\caption{The four effective degrees contributing at order $z^{-1}$ and their sector-valued contributions.}\label{tab:zminusone-degrees}
\begin{tabular}{c c c c}
\toprule
$(m,k)$ & $(d_0,d_1,d_2,d_3,d_u)$ & $d_F$ & contribution \\
\midrule
$(1,0)$ & $(\frac13,\frac23,1,\frac53,0)$ & $\frac{10}{3}$ & $14\,\1_\alpha$\\
$(0,1)$ & $(-\frac23,-\frac13,0,-\frac13,1)$ & $-\frac23$ & $\1_\alpha$\\
$(3,0)$ & $(1,2,3,5,0)$ & $10$ & $2520\,\1$\\
$(2,1)$ & $(0,1,2,3,1)$ & $6$ & $60\,\1$\\
\bottomrule
\end{tabular}
\end{table}

For each of these four degrees, Wang's localized class is the fundamental class of the displayed inertia sector.  Indeed, for $(1,0)$ and $(0,1)$ the hypersurface degree is nonintegral and no coordinate degree is a negative integer, while for $(3,0)$ and $(2,1)$ all coordinate degrees and the hypersurface degree are nonnegative integers.  Corollary~3.4 of \cite{Wang} therefore gives the asserted class in each case.

For $(1,0)$ the scalar coefficient is
\begin{equation}\label{eq:14}
\frac{(10/3)(7/3)(4/3)(1/3)}{(1/3)(2/3)(1)(5/3)(2/3)}=14,
\end{equation}
recovering the coefficient already present in the non-extended $I$-function of Oneto--Petracci \cite{OP}.  For $(0,1)$, all negative nonintegral products are empty and the positive degree $d_u=1$ contributes the unique denominator $z$, hence coefficient $1$.  For degree $(2,1)$, $d_0$ vanishes, so the $x_0$-coordinate contributes no denominator factor.  The remaining scalar coefficients are
\[
\frac{10!}{1!2!3!5!}=2520,
\qquad
\frac{6!}{1!2!3!}=60.
\]

Collecting the four contributions, the $z^{-1}$ part of the normalized $I$-function is
\begin{equation}\label{eq:I-leading}
I(q_1,q_2,z)
=
\1+
\frac{(2520q_1^3+60q_1^2q_2)\1+(14q_1+q_2)\1_\alpha}{z}
+O(z^{-2}).
\end{equation}
Lemma~\ref{lem:z-order} proves that \eqref{eq:I-leading} exhausts the nonnegative part of $zI-z\1$: every other nonzero effective summand begins at $z^{-2}$ or lower in $I$.  These four terms therefore determine the mirror transformation completely.

By the specialized restatement \cref{thm:wang-specialized} of Wang's published mirror theorem, the normalized $I$-function is identified with the $J$-function under the mirror transformation:
\[
I_{\cX_{10}}(q,0,z)=J_{\cX_{10}}(q,\mu(q,0,z),z),
\qquad
\mu(q,0,z)=\bigl[zI_{\cX_{10}}(q,0,z)-z\1\bigr]_+.
\]
Thus \eqref{eq:I-leading} gives the mirror input
\[
\mu=U\1+(14q_1+q_2)\1_\alpha,\qquad
U=2520q_1^3+60q_1^2q_2.
\]

The five coordinate characters in \eqref{eq:charge} sum to $(11,7)$, so subtracting the hypersurface character $(10,6)$ gives the extended anticanonical character $(1,1)$.  Since $u$ is nowhere zero on the stable locus, its tautological section trivializes $L_u$ on the quotient, and hence $L_{-K_{\cX_{10}}}\cong L_t\otimes L_u\cong L_t$.  A genuine stable-map class $d$ therefore has $d_u=0$.  If its degree in the extended lattice is $m\beta_1$, then its anticanonical degree is $e:=-K_{\cX_{10}}\cdot d=d_t=m/3$, since $\beta_1$ pairs with the anticanonical character by $1/3$.  Thus the ordinary Novikov monomial $q^e$ is represented by $q_1^m=q_1^{3e}$, which is the convention $q=q_1^3$ of \cite[Section 4]{OP}.  By contrast, an extended degree $\beta=m\beta_1+k\beta_2$ may have $k>0$ and need not be a stable-map degree.  Its anticanonical pairing is
\[
\beta(L_{(1,1)})=d_t+d_u=\frac{m+k}{3}.
\]
By \eqref{eq:period-degree}, an ordinary term has $t$-degree $e+n/3$, where $n$ is the number of insertions of $\1_\alpha$.  The anticanonical pairing of each of $\beta_1$ and $\beta_2$ is $1/3$, so the extended monomial $q_1^m q_2^k$ has degree $(m+k)/3$.  This gives a grading-compatible assignment of $t^{1/3}$ to both extended variables without identifying every extended degree with a stable-map degree.
In the identity sector, where $m=3\ell+2k$, this degree is $\ell+k$.  We choose the following specialization, with the MMLP parameter $s$ of degree zero:
\begin{equation}\label{eq:q-sub}
q_1=t^{1/3},
\qquad
q_2=s\,t^{1/3}.
\end{equation}
The variable $q_2$ records the additional $S$-extended direction and is not an ordinary stable-map Novikov variable.  This specialization sets $q_2=sq_1$ and sends the ordinary variable $q=q_1^3$ to $t$.  Substitution into the mirror input gives
\[
U\1+(14q_1+q_2)\1_\alpha=C(s)t\1+(14+s)t^{1/3}\1_\alpha.
\]
The twisted input after the mirror transformation is therefore $\lambda t^{1/3}\1_\alpha$, with twisted coordinate $\lambda$:
\begin{equation}\label{eq:lambda-s}
\lambda:=s+14=a-346.
\end{equation}

Under the specialization \eqref{eq:q-sub}, the unit input in the mirror transformation is
\[
U=2520q_1^3+60q_1^2q_2=C(s)t.
\]
Set $\tau=\lambda t^{1/3}\1_\alpha$.  Applying \cref{thm:wang-specialized} gives $I=J(\tau+U\1,z)$.  The normalized $J$-function satisfies the string equation
\[
J(\tau+U\1,z)=e^{U/z}J(\tau,z),
\]
so $J(\tau,z)=e^{-U/z}I$.  At $z=1$, this removes the unit input, multiplying the scalar identity-sector component by $e^{-C(s)t}$.

The scalar identity-sector component is extracted by pairing with the point class $\phi_{\mathrm{vol}}$.  A term of degree $d$ with $n$ insertions of $\tau$ contains the descendant contribution
\[
\frac{q^d}{n!}\left\langle\underbrace{\tau,\ldots,\tau}_{n},\frac{\phi_{\mathrm{vol}}}{z(z-\psi)}\right\rangle^{\cX_{10}}_{0,n+1,d},
\qquad
\frac{1}{z(z-\psi)}=\sum_{j\ge0}\psi^jz^{-j-2}.
\]
For $n$ insertions of the age-$2/3$ sector $\alpha$ and one final untwisted marking, the orbifold virtual dimension is
\[
\begin{aligned}
\operatorname{vdim}_{\C}\overline{\mathcal M}_{0,n+1}(\cX_{10},d)
&=(\dim_{\C}\cX_{10}-3)+(n+1)-n\operatorname{age}(\alpha)\\
&\qquad+c_1(T\cX_{10})\cdot d\\
&=-1+(n+1)-\frac{2n}{3}-K_X\cdot d\\
&=-K_X\cdot d+\frac n3=:r.
\end{aligned}
\]

The point class has complex codimension two, so the dimension constraint forces $j=r-2$.  The corresponding term therefore contributes $z^{-j-2}=z^{-r}$.
After specializing $q=t$ and $\tau=\lambda t^{1/3}\1_\alpha$ and setting $z=1$, we obtain the terms with $r\ge2$ in $G_{\cX_{10}}(\lambda;t)$.  By definition, its initial coefficients are $c_0(\lambda)=1$ and $c_1(\lambda)=0$.

The quantum period depends on the scalar identity-sector component.  A degree $m\beta_1+k\beta_2$ belongs to that sector precisely when
\[
m\equiv2k\pmod3.
\]
The parametrization $m=3\ell+2k$ makes the coordinate degrees especially simple.  When $\ell$ is nonnegative, they are the nonnegative integers
\[
\ell,\quad2\ell+k,\quad3\ell+2k,\quad5\ell+3k,\quad k.
\]
By Lemma~\ref{lem:effective-semigroup}, the class $\beta=(3\ell+2k)\beta_1+k\beta_2$ is effective for every $\ell,k\ge0$.  Moreover, the hypersurface degree is given by
\[
d_F=10\ell+6k.
\]
Corollary~3.4 of \cite{Wang} then identifies the localized Gysin contribution with the identity-sector fundamental class, since the coordinate degrees and $d_F$ are all nonnegative integers.  The resulting scalar coefficient is
\begin{equation}\label{eq:Clk}
C_{\ell,k}
=
\frac{(10\ell+6k)!}{\ell!\,k!\,(2\ell+k)!\,(3\ell+2k)!\,(5\ell+3k)!}
\end{equation}.

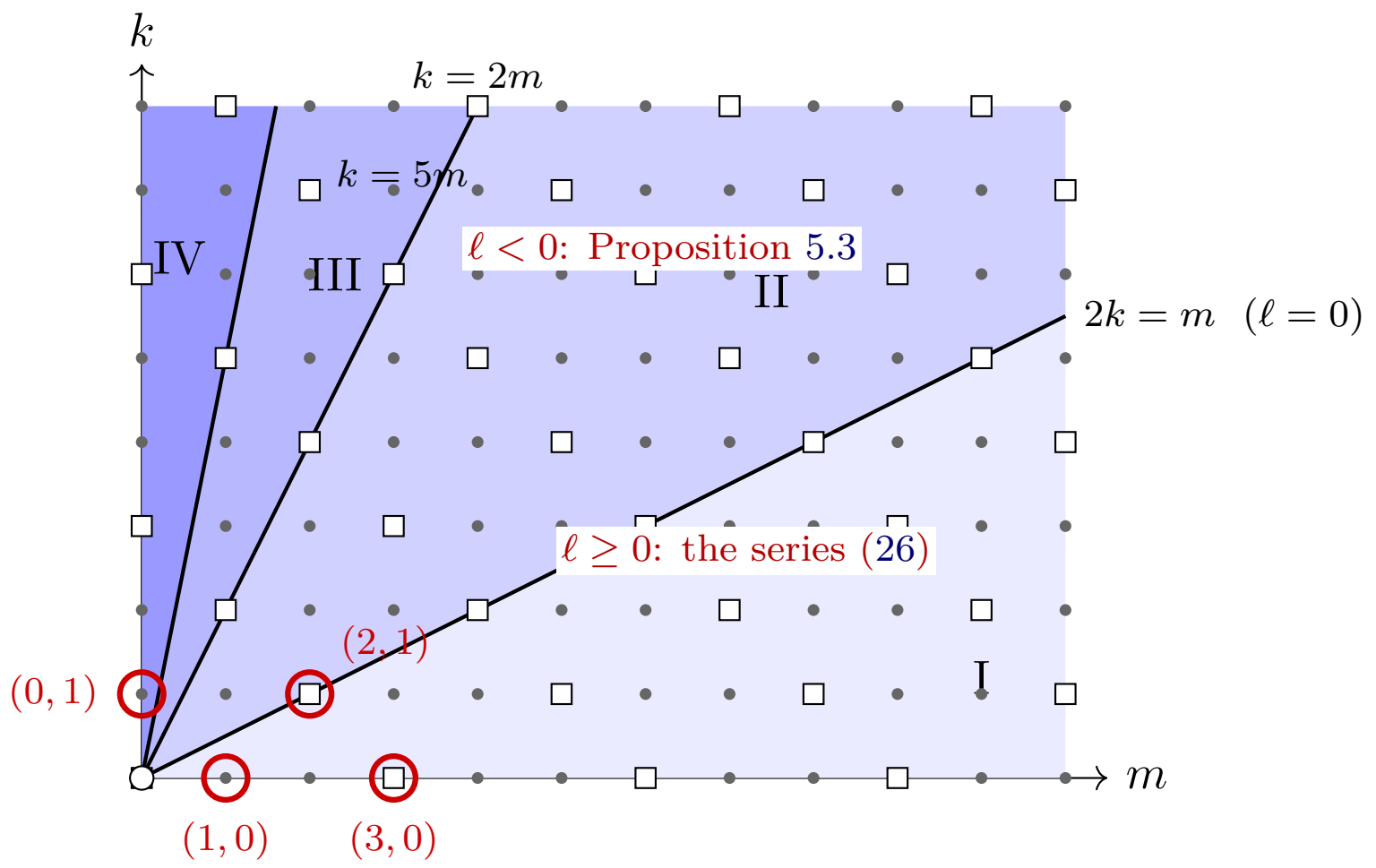
\begin{figure}[htbp]
\centering
\begin{tikzpicture}[scale=0.62]
  \draw[->] (0,0) -- (11.5,0) node[right,font=\small] {$m$};
  \draw[->] (0,0) -- (0,8.5) node[above,font=\small] {$k$};
  % regions shading
  \fill[blue!8]   (0,0) -- (11,0) -- (11,5.5) -- cycle;                 % I: 2k<=m
  \fill[blue!18]  (0,0) -- (11,5.5) -- (11,8) -- (4,8) -- cycle;        % II: m<2k, k<=2m
  \fill[blue!28]  (0,0) -- (4,8) -- (1.6,8) -- cycle;                   % III: 2m<k<=5m
  \fill[blue!40]  (0,0) -- (1.6,8) -- (0,8) -- cycle;                   % IV: 5m<k
  \draw[thick] (0,0) -- (11,5.5);   \node[font=\scriptsize,right] at (11,5.5) {$2k=m$ \ ($\ell=0$)};
  \draw[thick] (0,0) -- (4,8);      \node[font=\scriptsize,above] at (4,8) {$k=2m$};
  \draw[thick] (0,0) -- (1.6,8);    \node[font=\scriptsize,right] at (2.1,7.2) {$k=5m$};
  \node[font=\small] at (10,1.2) {I};
  \node[font=\small] at (7.5,5.8) {II};
  \node[font=\small] at (2.3,6) {III};
  \node[font=\small] at (0.45,6.2) {IV};
  % lattice points, identity sector m = 2k mod 3 marked as squares
  \foreach \m in {0,...,11}
    \foreach \k in {0,...,8}{
      \pgfmathtruncatemacro{\r}{mod(\m-2*\k+30,3)}
      \ifnum\r=0
        \draw[fill=white] (\m-0.12,\k-0.12) rectangle (\m+0.12,\k+0.12);
      \else
        \fill[black!60] (\m,\k) circle (0.07);
      \fi
    }
  % line ell = 0 : m = 2k
  \node[red!70!black,font=\scriptsize,fill=white,inner sep=1pt] at (6.2,6.3) {$\ell<0$: Proposition~\ref{prop:negative-ell}};
  \node[red!70!black,font=\scriptsize,fill=white,inner sep=1pt] at (7.2,2.7) {$\ell\ge0$: the series \eqref{eq:quantum-period}};
  % the four z^{-1} degrees
  \foreach \p in {(1,0),(0,1),(3,0),(2,1)} \draw[red!80!black,very thick] \p circle (0.26);
  \node[red!80!black,font=\scriptsize,below] at (1,-0.3) {$(1,0)$};
  \node[red!80!black,font=\scriptsize,left] at (-0.3,1) {$(0,1)$};
  \node[red!80!black,font=\scriptsize,below] at (3,-0.3) {$(3,0)$};
  \node[red!80!black,font=\scriptsize,above right] at (2.15,1.15) {$(2,1)$};
  \fill[white] (0,0) circle (0.14); \draw (0,0) circle (0.14);
\end{tikzpicture}
\caption{Effective degrees $(m,k)$ of the $S$-extended presentation. The walls $2k=m$, $k=2m$, and $k=5m$ cut the quadrant into the regions I--IV of Table~\ref{tab:wall-regions}. Squares mark the identity sector $m\equiv2k\pmod3$, and dots mark the $\mu_3$ sectors. The four circled degrees are the only ones whose summand begins at order $z^{-1}$ (Lemma~\ref{lem:z-order}). With $m=3\ell+2k$, the wall $2k=m$ is the line $\ell=0$. The identity-sector points of Region~I give the scalar series \eqref{eq:quantum-period}, while those of Regions II--IV have zero scalar projection by Proposition~\ref{prop:negative-ell}.}
\label{fig:degrees}
\end{figure}

Figure~\ref{fig:degrees} summarizes the wall regions, the identity-sector lattice, and the boundary between nonnegative and negative $\ell$.

\FloatBarrier

Effectivity does not force $\ell\ge0$, so degrees with $\ell<0$ could a priori contribute additional identity-sector terms.  The following proposition proves that the displayed scalar contribution is complete.  Set $D:=\cX_{10}\cap\{x_0=0\}$.

\begin{proposition}[Completeness of the scalar identity-sector contribution]\label{prop:negative-ell}
Suppose that $m$ has the form $3\ell+2k$ with $\ell<0$.  Then the localized Gysin class in Wang's $(m,k)$-summand has zero projection to the scalar $H^0$ of the identity sector.  Consequently, the scalar identity-sector contribution is obtained exactly from $\ell,k\ge0$.
\end{proposition}

\begin{proof}
Here $g_\beta$ is the identity element, so all coordinate degrees and the hypersurface degree $d_F=10\ell+6k$ are integers.
Write $Y_\beta^{ss}$ for the corresponding hypersurface fixed locus.
Since $d_0=\ell<0$, Wang's fixed-locus construction sets $x_0=0$ on $Z_\beta^{ss}$ \cite[Section 3, (3.1), and Definition 3.3]{Wang}.
If $d_F<0$, then $d_F$ is a negative integer.  Wang's Lemma~3.2 gives the vanishing of the hypersurface equation on the distinguished fixed locus $Z_\beta^{ss}$, so the hypersurface fixed locus is all of $Z_\beta^{ss}$.  Thus
\[
Z_\beta^{ss}=Y_\beta^{ss}\hookrightarrow\{x_0=0\}\cap\cX_{10}=D.
\]
If $d_F\ge0$, Wang imposes the defining section, and
\[
Y_\beta^{ss}\subseteq Z_\beta^{ss}\cap\{\widetilde F=0\}\subseteq D.
\]
In either case, the source of the localized Gysin class maps through $D$, so its pushforward factors through the inclusion
\[
j:D\hookrightarrow\cX_{10}.
\]
After expanding the hypergeometric factors in $z^{-1}$, their coefficients are polynomials in the divisor classes $D_\rho$, so have nonnegative cohomological degree.

\Needspace{4\baselineskip}
The projection formula gives $h\cup j_*\gamma=j_*(j^*h\cup\gamma)$ for each such coefficient $h$.  Thus multiplication preserves support on $D$ and cannot produce an $H^0$-component.  The substack $D$ has positive codimension because the defining polynomial of a general $X_{10}$ does not vanish identically on $x_0=0$, as it contains the nonzero term $x_3^2$.  Thus the full summand has zero scalar projection, although the localized Gysin class itself need not vanish.
\end{proof}

Combining \eqref{eq:I-leading}--\eqref{eq:Clk} gives
\begin{equation}\label{eq:quantum-period}
\boxed{
G_{\cX_{10}}(\lambda;t)
=
 e^{-C(s)t}
\sum_{\ell,k\ge0}
\frac{(10\ell+6k)!}{\ell!\,k!\,(2\ell+k)!\,(3\ell+2k)!\,(5\ell+3k)!}
 s^k t^{\ell+k},
}
\end{equation}
where $s=\lambda-14$ as in \eqref{eq:lambda-s}.

It remains to compare \eqref{eq:quantum-period} with the classical period of $f_s$.  The proof identifies its factorial coefficients as constant terms and then applies the multinomial theorem.

\begin{proof}[Proof of Theorem~\ref{thm:conjB}]
For any $\ell,k\ge0$, the product of the previously defined factors is
\[
A^\ell B^k=x^{-3\ell-2k}y^{5\ell+3k}H^{10\ell+6k}.
\]
Its constant term is obtained by choosing from $H^{10\ell+6k}=(1+x+y^{-1})^{10\ell+6k}$ precisely $3\ell+2k$ copies of $x$, $5\ell+3k$ copies of $y^{-1}$, and $2\ell+k$ copies of $1$.  These multiplicities sum to $10\ell+6k$, so the multinomial theorem gives
\begin{equation}\label{eq:CTAB}
\CT(A^\ell B^k)
=\frac{(10\ell+6k)!}{(3\ell+2k)!(5\ell+3k)!(2\ell+k)!}.
\end{equation}
The summand in \eqref{eq:quantum-period} is therefore $\CT(A^\ell B^k)/(\ell!k!)$.  Expanding the exponential gives
\begin{equation}\label{eq:classical-coeff}
G_{\cX_{10}}(\lambda;t)
=\sum_{\ell,k,r\ge0}
\frac{s^k\bigl(-C(s)\bigr)^r}{\ell!k!r!}
\CT(A^\ell B^k)t^{\ell+k+r}.
\end{equation}
After regularization, the coefficient of $t^n$ is
\[
n!\sum_{\ell+k+r=n}
\frac{s^k\bigl(-C(s)\bigr)^r}{\ell!k!r!}\CT(A^\ell B^k),
\]
which is exactly $\CT(f_s^n)$ by the multinomial theorem.  This coefficient identity holds for every $n\ge0$ and every $s\in\A^1$, so the regularized quantum period and the classical period agree as full formal power series for every member of the normalized family.  It follows that
\[
\widehat G_{\cX_{10}}(\lambda;t)=\pi_{f_s}(t).
\]
Equation~\eqref{eq:lambda-s} therefore gives the desired identity $\widehat G_{\cX_{10}}(\lambda;t)=\pi_{f_a}(t)$.
By Lemma~\ref{lem:parameter-dimensions}, both $L_P^T$ and $H^{\mathrm{ts}}_{\cX_{10}}$ are one-dimensional.  Since $a\mapsto a-346$ has nonzero linear part, the map
\[
\Phi(f_a)=(a-346)\1_\alpha
\]
is an affine isomorphism.  This proves Theorem~\ref{thm:conjB}.
\end{proof}

\begin{remark}[Normalization checks]\label{rem:normalization-checks}
At the Oneto--Petracci point, $a=360$, so $s=0$ and $\lambda=14$.  Specializing \eqref{eq:quantum-period} at $s=0$ and regularizing gives
\[
1+32448360t^2+515050578000t^3+\cdots,
\]
which is the period recorded in \cite[entry 11]{OP}.  At $\lambda=0$, equation~\eqref{eq:lambda-s} gives $s=-14$, hence $a=346$ and
\[
f_{346}=A-14B-1680.
\]
\end{remark}

\section*{Further directions}

The Corti--Heuberger classification contains eight Fano-index-one families $X_{1,d}$, each with a unique $\frac13(1,1)$ singularity.  Thus the twisted parameter space relevant to Conjecture~B remains one-dimensional throughout this cascade \cite[Corollary 8 and Table 2]{CortiHeuberger}.  The neighboring family $X_{1,4/3}$ is a natural next test, and a general member is a complete intersection of two quartics in $\Pj(1,1,2,2,3)$.  This setting differs from the nonabelian weighted-Grassmannian construction used by Sultani--Webb for $X_{1,7/3}$ \cite{SultaniWebb}.

\section*{Funding}

This research was supported by FAPERJ (Fundação Carlos Chagas Filho de Amparo à Pesquisa do Estado do Rio de Janeiro) through the Programa Bolsa Mestrado e Doutorado Nota 10 – 2025 under Chamada FAPERJ Nº 01/2025 (Processo SEI-260003/001060/2025).

\section*{Declaration of generative AI and AI-assisted technologies in the writing process}

During the preparation of this manuscript, the author used ChatGPT (OpenAI) for language editing, organization, \LaTeX{} formatting, bibliographic research, mathematical exposition, and exploratory symbolic and numerical checks.  OpenAI Codex was used in the final presentation pass for proofreading, bibliography verification, and PDF layout inspection.

\AtNextBibliography{\small}
\printbibliography

\end{document}